\documentclass[12pt]{amsart}
\usepackage[top=1in, bottom=1in, left=1in, right=1in]{geometry}
\usepackage{dsfont}
\usepackage{amsmath} 
\usepackage{amsthm}
\usepackage{amssymb} 
\usepackage{mathtools}
\usepackage{times}
\usepackage[utf8]{inputenc}
\usepackage{indentfirst}
\usepackage{graphicx} 
\usepackage{color}
\usepackage[colorlinks=true]{hyperref}
\hypersetup{colorlinks=true, citecolor=green, linkcolor=green, filecolor=magenta, urlcolor=cyan}
\usepackage[shortlabels]{enumitem} 
\usepackage{mdframed}
\usepackage{tcolorbox}
\newtcolorbox{mybbox}{colback=blue!5!white,
colframe=blue!75!black}
\newtcolorbox{mygbox}{colback=green!5!white,
colframe=green!75!black}
\newtcolorbox{myrbox}{colback=red!5!white,
colframe=red!75!black}

\newtheorem{theorem}{Theorem}[section]
\newtheorem*{thm*}{Theorem FI (Friedlander-Iwaniec)}

\newtheorem{corollary}[theorem]{Corollary}

\newtheorem{lemma}[theorem]{Lemma}

\newtheorem{proposition}[theorem]{Proposition}

\theoremstyle{remark}

\numberwithin{equation}{section}

\newcommand{\N}{\mathbb{N}}

\newcommand{\R}{\mathbb{R}}

\newcommand{\CD}{\mathcal{D}}

\DeclareMathOperator{\supp}{supp}

\newcommand{\nn}{\nonumber \\}

\newcommand{\pthbgg}[1]{\bigg( #1 \bigg)}

\newcommand{\abs}[1]{\lvert #1 \rvert}

\newcommand{\absBgg}[1]{\Bigg\lvert #1 \Bigg\rvert}

\newcommand{\flbgg}[1]{\bigg\lfloor #1 \bigg\rfloor}

\newcommand{\bs}\boldsymbol{}
\renewcommand{\geq}{\geqslant}
\renewcommand{\leq}{\leqslant}

\renewcommand{\bar}[1]{\overline{#1}}

\newcommand{\eps}{\varepsilon}
\renewcommand{\Re}{{\rm Re}}

\renewcommand{\mod}[1]{\,({\rm mod}\,#1)}

\UseRawInputEncoding
\begin{document}

\title{Primes in arithmetic progressions under the presence of Landau-Siegel zeroes}

\author{Stelios Sachpazis}

\address{Charles University, Faculty of Mathematics and Physics, Department of Algebra, Sokolovsk\'{a} 49/83, 186 75 Praha 8, Czech Republic}
\email{stylianos.sachpazis@matfyz.cuni.cz}


\date{\today}

\maketitle

\begin{abstract}
Let $x\geq 2$ and assume that $a$ and $q$ are coprime positive integers. As usual, $\psi(x;q,a)\!\vcentcolon=\sum_{n\leq x,n\equiv a\mod{q}}\Lambda(n)$, where $\Lambda$ is the von Mangoldt function. In $2003$, Friedlander and Iwaniec assumed the existence of exceptional characters corresponding to ``extreme" Landau-Siegel zeroes and established a meaningful asymptotic formula for $\psi(x;q,a)$ beyond the square-root barrier of the Generalized Riemann Hypothesis. In particular, their asymptotic yields non-trivial information for moduli $q\leq x^{1/2+1/231}$. In this paper, we considerably relax the ”extremity” of the Landau-Siegel zero required in the work of Friedlander and Iwaniec and obtain a conditional asymptotic formula for $\psi(x;q,a)$ in a slightly wider range of $q$.
\end{abstract}

\section{Introduction}

Let $x\geq 2$ and $\Lambda$ be the von Mangoldt function. Assume that $a$ and $q$ are two coprime positive integers and let, as usual,
\begin{equation*}
\psi(x;q,a)\vcentcolon=\sum_{\substack{n\leq x\\n\equiv a\mod{q}}}\Lambda(n).
\end{equation*}

The Siegel-Walfisz theorem, frequently referred to as the prime number theorem in arithmetic progressions, states that given $A>0$, there exists some positive constant $c_A$ depending on $A$ such that
\begin{equation}\label{sw}
\psi(x;q,a)=\frac{x}{\phi(q)}+O\Big(xe^{-c_A\sqrt{\log x}}\Big)
\end{equation}
for all $q\leq (\log x)^A$. Here and hereafter, the Greek letter $\phi$ denotes the Euler totient function.

The overall quality of the Siegel-Walfisz theorem is tied to our knowledge on the location of the zeroes of the Dirichlet $L$-functions. To be more precise, the wider zero-free regions one can establish for the Dirichlet $L$-functions, the better error term and range of $q$ they can obtain. For example, if we assume the Generalized Riemann Hypothesis (abbreviated as GRH for short), which claims that all the zeroes of the Dirichlet $L$-functions lie on the vertical line $\sigma=1/2$, then the error term and the uniformity (range of $q$ for which the derived asymptotics of $\psi(x;q,a)$ are non-trivial) in the Siegel-Walfisz theorem improve substantially. In fact, GRH yields
\begin{equation}\label{GRH}
\psi(x;q,a)=\frac{x}{\phi(q)}+O(\sqrt{x}(\log x)^2)
\end{equation}
for every $q\leq x$.

Note that (\ref{GRH}) constitutes a meaningful asymptotic formula as long as its anticipated main term $x/\phi(q)$ dominates the error term. This is clearly the case for $q\leq x^{1/2-\eps}$, but for $q$ exceeding $\sqrt{x}$, the expression $\sqrt{x}(\log x)^2$ in the big-Oh term ``overpowers" $x/\phi(q)$. 

Summarizing, if one decides to approach the prime number theorem in arithmetic progressions conditionally, and appoints GRH as their main conjectural assumption, then they can considerably decrease the unconditional error term in (\ref{sw}) and derive a meaningful asymptotic formula for all $q$ residing in the much wider range $[1,x^{1/2-\eps}]$.

Upon witnessing these strong refinements that GRH implies over the Siegel-Walfisz theorem, one might naturally wonder what more could be said for primes in arithmetic progressions under other unresolved conjectures concerning the zeroes of the Dirichlet $L$-functions. In this direction, we could ask whether a non-trivial asymptotic formula for $\psi(x;q,a)$ can hold for moduli $q$ beyond the $\sqrt{x}$-limitation of GRH under the presence of Landau-Siegel zeroes.

A \textit{Landau-Siegel zero}, often just called a \textit{Siegel zero}, is a real number $\beta$ associated to a real primitive Dirichlet character $\chi$ modulo $D\geq 2$ such that $L(\beta,\chi)=0$ and 
$$\beta=1-\frac{1}{\eta\log D}$$
for some $\eta\geq 2$. The number $\eta$ is called the \textit{quality} of the Landau-Siegel zero $\beta$ and from Siegel's theorem we have the ineffective bound
\begin{equation}\label{siegel}
    \eta\ll_{\eps}D^{\eps}
\end{equation}
for any $\eps>0$.

Even though it is widely believed that Landau-Siegel zeroes are unlikely to exist, their improbable presence can have intriguing implications. For instance, from a 1983 result of Heath-Brown \cite{heath}, it follows that if there exist infinitely many Landau-Siegel zeroes $\beta_j=1-1/(\eta_j\log D_j)$ such that $\eta_j\to +\infty$, then the Twin Prime Conjecture must be true.

Let us now return to the very reason for which Landau-Siegel zeroes entered our discussion in the first place. We would like to asymptotically evaluate $\psi(x;q,a)$ for $q>\sqrt{x}$ under the existence of exceptional Dirichlet characters corresponding to Landau-Siegel zeroes. The first result of the literature on this problem came from Friedlander and Iwaniec who proved the following in $2003$.

\begin{thm*}
Let $a$ and $q$ be two coprime positive integers and assume that $\chi$ is a real primitive Dirichlet character modulo $D\geq 2$. Set also $r=554,401$ for simplicity. Then for $x\geq D^r$ and $q\leq x^{1/2+1/231}=x^{233/462}$, we have
\begin{equation*}
\psi(x;q,a)=\frac{\psi(x)}{\phi(q)}\bigg(1-\mathds{1}_{D\mid q}\chi(a)+O\Big(L(1,\chi)(\log x)^{r^r}\Big)\bigg),
\end{equation*}
where $\psi(x)\vcentcolon=\sum_{n\leq x}\Lambda(n)$ and $L(1,\chi)=\sum_{n\geq 1}\chi(n)n^{-1}$.
\end{thm*}

Note that if $x$ is a high power of $D$ in Theorem FI, then the theorem of Friedlander and Iwaniec furnishes a non-trivial asymptotic formula when $L(1,\chi)$ is very small with respect to $D$, e.g. when $L(1,\chi)\leq (\log D)^{-r^r-1}$. If we now assume that the character $\chi$ has a Landau-Siegel zero of quality $\eta$, then according to the work of Friedlander and Iwaniec from \cite{fin}, we have that 
$$L(1,\chi)\ll\frac{\log D}{\eta\log\log D},$$
which is the best current bound on $L(1,\chi)$ in terms of $\eta$.
This means that if we translated Theorem FI into terms of $\eta$, the theorem would achieve its anticipated objective for $\eta\geq (\log D)^{r^r+2}$, say. Simply put, in the language of Landau-Siegel zeroes, the result of Friedlander and Iwaniec surpasses the $\sqrt{x}$-barrier of GRH and returns a sensible asymptotic formula for $\psi(x;q,a)$ with $q$ as large as $x^{233/462}$ when the character $\chi$ is accompanied by an ``extreme" Landau-Siegel zero. The word ``extreme" in this context means that the quality $\eta$ of the respective Landau-Siegel zero should be growing to infinity at a faster rate than a colossal power of $\log D$ (recall that $r$ above was about half a million). In a recent preprint, Wright \cite{wr} claims that he can relax the extremity of the Landau-Siegel zero required in the theorem of Friedlander and Iwaniec and obtain a non-trivial asymptotic formula for $\psi(x;q,a)$ under the assumption $L(1,\chi)=o((\log x)^{-7})$ as $x\to\infty$. Hence, if we assume that $x$ is a sufficiently large power of $D$, then the quality $\eta$ of the respective Landau-Siegel zero in Wright's work has to satisfy an inequality such as $\eta\geq (\log D)^A$ for some $A>8$.

The size conditions on $\eta$, under which Landau-Siegel-zero versions of Theorem FI and Wright's work deliver a non-trivial asymptotic formula for $\psi(x;q,a)$, seemed to admit improvement and this motivated the making of the current paper. More precisely, in this work, we show that a sensible asymptotic formula of $\psi(x;q,a)$ can be attained under no assumptions on the size of $\eta$ with respect to $D$.

\begin{theorem}\label{maint}
Let $x\geq 2,\eps\in(0,1/100)$, and consider positive integers $a,q$ and $D\geq 2$ such that $(a,q)=1,q\leq x^{58/115-\eps},$ and $x=D^V$ for some $V\geq 200/\eps$. Let also $\chi$ be a quadratic primitive character $\mod{D}$ and assume that $\eta_0=\eta_0(\eps)>0$ is a sufficiently large real number in terms of $\eps$. If for some $\eta\geq \eta_0$, there exists a real number $\beta=1-1/(\eta\log D)$ such that $L(\beta,\chi)=0$, then
\begin{align}\label{thmas}
\psi(x;q,a)=\frac{\psi(x)}{\phi(q)}\bigg\{1-\mathds{1}_{D\mid q}\chi(a)+O_{\eps}\bigg(\frac{V^{16}}{\eta}+\exp\Big(-C_{\eps}\sqrt{V\log\eta}\Big)\bigg)\bigg\},
\end{align}
where $C_{\eps}$ is a positive constant that depends on $\eps$.
\end{theorem}

Theorem \ref{maint} has the following immediate corollary.

\begin{corollary}\label{cor}
Under the same considerations and notation as in Theorem \ref{maint}, we have the following:
\begin{enumerate}[(a)]
    \item There exists a positive constant $C'_{\eps}$ depending on $\eps$ such that
\begin{align*}
\psi(x;q,a)=\frac{\psi(x)}{\phi(q)}\bigg\{1-\mathds{1}_{D\mid q}\chi(a)+O_{\eps}\bigg(\!\exp\Big(-C_{\eps}'\sqrt{\log\eta}\Big)\bigg)\bigg\},
\end{align*}
     for every $x\leq D^{C_{\eps}^{-2}\log\eta}$.
     \item For every $\delta\in(0,1)$, we have
\begin{eqnarray*} 
\psi(x;q,a)=\frac{\psi(x)}{\phi(q)}\bigg\{1-\mathds{1}_{D\mid q}\chi(a)+O_{\eps}\bigg(\frac{1}{\eta^{1-\delta}}\bigg)\bigg\}
\end{eqnarray*}
      for all $x \in (D^{{C_{\eps}}^{-2}\log{\eta}},D^{\eta^{\delta/16}}]$.
\end{enumerate}
\end{corollary}

As Friedlander and Iwaniec claimed in \cite{fip}, Theorem FI could be made stronger by pre-sieving. It turns out that their claim was spot on as we were able to mollify the exceptionality of the characters required in their theorem by the introduction of a preliminary sieve.

The error term in Theorem \ref{maint} is similar to the error terms obtained by Matom\"{a}ki and Merikoski in \cite{mm} and the bound that Jaskari and the author established for the Chowla-type $k$-point correlations of the Liouville function in \cite{js}. The reason for this lies in the common use of Lemma \ref{joka} that Matom\"{a}ki and Merikoski proved in \cite{mm} and in the sifting arguments that are similar to the ones that were first developed by Matom\"{a}ki and Merikoski in their precited work.

As was mentioned above, we would like to underline that this article's main scope was to loosen the nature of the Landau-Siegel zeroes needed in the result of Friedlander and Iwaniec and in the work of Wright. We did not attempt to broaden the range of $q$ to $[1,x^{\theta}]$ for some $\theta\geq58/115=1/2+1/230$, but we might consider looking into this possibility in the future. In his aforementioned preprint, Wright claims the wider range $q\leq x^{1/2+1/118-\eps}=x^{30/59-\eps}$. We believe that one could improve on this and establish Theorem \ref{maint} with an even wider range on $q$. A seemingly promising strategy for achieving this goal could be based off Heath-Brown's treatment on the ternary divisor function in arithmetic progressions from \cite{hbd3}. We presume that an appropriate adjustment of Heath-Brown's work -- if feasible -- could increase the exponent of $x$ in Lemma \ref{compa} to $1/2+1/82-\eps=21/41-\eps$. This would in turn enlarge our theorem's range to $q\in[1,x^{21/41-\eps}]$ for general moduli $q$.

\subsection*{Notation}
In this text, we use the classical notation $e(t)=e^{2\pi it}$ for $t\in\R$. For $a,b\in\N$, every occurrence of $(a,b)$ denotes their greatest common divisor. Throughout the paper, given a positive integer $n>1$, the notations $P^-(n)$ and $P^+(n)$ stand for the smallest and largest prime factor of $n$, respectively. If $n=1$, we conventionally have that $P^-(1)=+\infty$ and that $P^+(1)=0$. There are also instances where we encounter or mention the divisor function $\tau$, which is the arithmetic function defined as $\tau(n)=\sum_{d\mid n}1$ for all $n\in\N$. We also refer to the ternary divisor function $\tau_3$ which is the arithmetic function given by $\tau_3(n)=\sum_{abc=n}1$ for all $n\in\N$. Furthermore, given a non-principal Dirichlet character $\chi$, we use the notation $L(\cdot,\chi)$ for its Dirichlet $L$-function defined as $L(s,\chi)=\sum_{n\geq 1}\chi(n)n^{-s}$ in the open half-plane $\Re(s)>0$. Lastly, for arithmetic functions $f$ and $g$, the notation $f\ast g$ stands for their Dirichlet convolution, which is the arithmetic function given as $(f\ast g)(n)=\sum_{ab=n}f(a)g(b)$ for all $n\in\N$.

\subsection*{Acknowledgements}
The author would like to thank Kaisa Matom\"{a}ki for bringing this problem to his attention and for their discussions on the topic. During the making of this article, the author received support from the Finnish Research Council grant nos. 333707, 346307 and 364214. He also acknowledges support from the two Charles University programmes PRIMUS/25/SCI/008 and PRIMUS/25/SCI/0017. The biggest part of the present work was undertaken at the Department of Mathematics and Statistics of the University of Turku. The author would like to deeply thank the department for the excellent working environment during his time there.

\section{Heuristics, the initial reduction, and proof description}

\subsection{The opening heuristics} We start the current section by providing an intuitive consequence of a Landau-Siegel zero's presence. It will be a non-rigorous conclusion which serves as the basis of all the conditional works on primes that build upon the existence of Landau-Siegel zeroes. 

Let $\chi$ be a quadratic primitive Dirichlet character modulo $D\geq 2$ such that there exists a Landau-Siegel zero $\beta=1-1/(\eta\log D)$, where we think of $\eta$ as being large. By the continuity of $L(\cdot,\chi)$, and since $\eta$ is considered large, we expect $L(1,\chi)$ to be ``close" to $L(\beta,\chi)=0$. Therefore, $L(1,\chi)$ is forced to be small when $\chi$ is exceptional. However,
$$L(1,\chi)^{-1}=\prod_p\bigg(1-\frac{\chi(p)}{p}\bigg),$$
and so the Euler product of the right-hand side has to be large. In light of Mertens' theorem, this means that $\chi(p)$ should prefer the value $-1$ on a considerable amount of the large primes. In other words, the character $\chi$ ``pretends" to behave like the M\"{o}bius function $\mu$ on many large primes. It then seems safe to expect that $\Lambda=\mu\ast\log\approx\chi\ast\log$ on those large primes. Taking this heuristically driven conclusion a step further, we claim that under the presence of Landau-Siegel zeroes, the arithmetic function $$\lambda'\vcentcolon=\chi\ast\log$$ 
can be viewed as a sensible approximation of $\Lambda$ on integers with large prime factors. This claim is the foundation of treatments that study the primes under the unlikely existence of Landau-Siegel zeroes. 

\subsection{The initial reduction for the proof of Theorem \ref{maint}}
According to the concluding assertion of the previous subsection, if $\chi$ is an exceptional Dirichlet character possessing a Landau-Siegel zero of high quality, then we believe that $\chi\ast\log=\lambda'\approx\Lambda$ on integers with appropriately large prime factors. Consequently, the existence of Landau-Siegel zeroes suggests that instead of studying $\sum_{n\leq x,n\equiv a\mod{q}}\Lambda(n)$ directly, we can first drop the negligible contribution coming from the small primes and their powers, and then attempt to establish an asymptotic formula for the resulting sums of $\Lambda$ via the sums
$$\sum_{\substack{n\leq x\\n\equiv a\mod{q}\\P^-(n)>z}}\lambda'(n),$$
where $z$ is some suitable sifting parameter that will be eventually chosen in such a way so that it is at most as large as a small power of $x$. The realization and execution of this idea toward a proof of Theorem \ref{maint} necessitate the attainment of a sound comparison between
\begin{align}\label{2sums}
\sum_{\substack{n\leq x\\n\equiv a\mod{q}\\P^-(n)>z}}\Lambda(n)\quad\text{ and}\quad\!\!\sum_{\substack{n\leq x\\n\equiv a\mod{q}\\P^-(n)>z}}\lambda'(n).
\end{align}

In order to establish a relation involving the sums (\ref{2sums}), we are in search of a tangible connection between $\lambda'$ and $\Lambda$. Such a connection can be derived from the definition of $\lambda'$. Indeed, upon setting 
$$\lambda\vcentcolon=1\ast\chi,$$ 
we see that
$$\lambda'=\chi\ast\log=\chi\ast(1\ast\Lambda)=(1\ast\chi)\ast\Lambda=\lambda\ast\Lambda,$$ 
which implies that
\begin{equation}\label{lL}
\lambda'(n)=\Lambda(n)+\sum_{\substack{k\ell=n\\k>1}}\lambda(k)\Lambda(\ell).
\end{equation}

We now apply (\ref{lL}) to concretely compare the two sums in (\ref{2sums}). In particular, (\ref{lL}) yields
\begin{align}\label{sp}
\sum_{\substack{n\leq x\\n\equiv a\mod{q}\\P^-(n)>z}}\lambda'(n)=\sum_{\substack{n\leq x\\n\equiv a\mod{q}\\P^-(n)>z}}\Lambda(n)+\sum_{\substack{k\ell\leq x,\,k>z\\k\ell\equiv a\mod{q}\\P^-(k\ell)>z}}\lambda(k)\Lambda(\ell).
\end{align}

This should be the starting point for the proof of Theorem \ref{maint}, but we would like it to involve the original object of study, that is, the quantity $\psi(x;q,a)$. We thus replace the sifted sum of $\Lambda$ over the arithmetic progression $a\mod{q}$ with $\psi(x;q,a)$ by using the following elementary estimate:
\begin{align}\label{bes}
\sum_{\substack{n\leq x\\P^-(n)\leq z}}\Lambda(n)\leq \sum_{p\leq z}\log p\sum_{\nu:p^{\nu}\leq x}1=\sum_{p\leq z}\log p\flbgg{\frac{\log x}{\log p}}\leq \pi(z)\log x\ll \frac{z\log x}{\log z}.
\end{align}
So, relation (\ref{sp}) becomes
\begin{align}\label{sprv}
\sum_{\substack{n\leq x\\n\equiv a\mod{q}\\P^-(n)>z}}\lambda'(n)=\psi(x;q,a)+\sum_{\substack{k\ell\leq x,\,k>z\\k\ell\equiv a\mod{q}\\P^-(k\ell)>z}}\lambda(k)\Lambda(\ell)+O\bigg(\frac{z\log x}{\log z}\bigg).
\end{align}

This is our starting point revisited, and we now examine how it can lead to (\ref{thmas}). With this goal in mind, we first notice that locating the ``main term" $\psi(x)/\phi(q)$ in (\ref{sprv}) appears to be an elusive task. To this end, let us force it into existence. Upon a second application of (\ref{lL}), the expression of interest, namely $\psi(x)/\phi(q)$, is found in
\begin{align}\label{cmt}
\begin{split}
\frac{1}{\phi(q)}\sum_{\substack{n\leq x\\(n,q)=1\\P^-(n)>z}}\lambda'(n)&=\frac{1}{\phi(q)}\sum_{\substack{n\leq x\\(n,q)=1\\P^-(n)>z}}\Lambda(n)+\frac{1}{\phi(q)}\sum_{\substack{k\ell\leq x,\,k>z\\(k\ell,q)=1\\P^-(k\ell)>z}}\lambda(k)\Lambda(\ell)\\
&=\frac{\psi(x)}{\phi(q)}+\frac{1}{\phi(q)}\sum_{\substack{k\ell\leq x,\,k>z\\(k\ell,q)=1\\P^-(k\ell)>z}}\lambda(k)\Lambda(\ell)+O\bigg(\frac{z\log x}{\phi(q)\log z}\bigg),
\end{split}
\end{align}
where we made use of the basic estimate (\ref{bes}) again. 

Multiplying both sides of (\ref{cmt}) by $1-\mathds{1}_{D\mid q}\chi(a)$, and subtracting from (\ref{sprv}), we deduce that
\begin{align}\label{tocon}
\begin{split}
\psi(x;q,a)-\frac{\psi(x)}{\phi(q)}(1-\mathds{1}_{D\mid q}\chi(a))=&\sum_{\substack{n\leq x\\n\equiv a\mod{q}\\P^-(n)>z}}\lambda'(n)-\frac{1-\mathds{1}_{D\mid q}\chi(a)}{\phi(q)}\sum_{\substack{n\leq x\\(n,q)=1\\P^-(n)>z}}\lambda'(n)\\
&+\frac{1-\mathds{1}_{D\mid q}\chi(a)}{\phi(q)}\sum_{\substack{k\ell\leq x,\,k>z\\(k\ell,q)=1\\P^-(k\ell)>z}}\lambda(k)\Lambda(\ell)\\
&-\sum_{\substack{k\ell\leq x,\,k>z\\k\ell\equiv a\mod{q}\\P^-(k\ell)>z}}\lambda(k)\Lambda(\ell)+O\bigg(\frac{z\log x}{\phi(q)\log z}\bigg).
\end{split}
\end{align}

Consequently, it becomes clear that the proof of Theorem \ref{maint} reduces to establishing bounds on the difference
\begin{align}\label{Delta}
\Delta\vcentcolon=\sum_{\substack{n\leq x\\n\equiv a\mod{q}\\P^-(n)>z}}\lambda'(n)-\frac{1-\mathds{1}_{D\mid q}\chi(a)}{\phi(q)}\sum_{\substack{n\leq x\\(n,q)=1\\P^-(n)>z}}\lambda'(n),
\end{align}
and the double sums
\begin{align}
\label{S1}S_1\vcentcolon=&\sum_{\substack{k\ell\leq x,\,k>z\\(k\ell,q)=1\\P^-(k\ell)>z}}\lambda(k)\Lambda(\ell),\\
\label{S2}S_2\vcentcolon=&\sum_{\substack{k\ell\leq x,\,k>z\\k\ell\equiv a\mod{q}\\P^-(k\ell)>z}}\lambda(k)\Lambda(\ell).
\end{align}
The bounds will depend on the Landau-Siegel zero and yield the error term stated in Theorem \ref{maint}.

\subsection{A brief description of the proof of Theorem \ref{maint}}
Let us outline in a few words how we are going to bound $\Delta$ and the sums $S_1,S_2$ in the presence of Landau-Siegel zeroes.

Concerning $\Delta$, we first detect and remove the roughness condition $P^-(\cdot)>z$ by implementing classical sieve weights. We then apply Lemma \ref{l's} and Proposition \ref{slaps} to the respective occurring sums of $\lambda'$ and $\lambda$ over arithmetic progressions. We proceed to make use of the fundamental lemma of sieve theory along with a variant of it containing logarithmic weights (for the latter, see Lemma \ref{FL}(b)), and then, we finish off the estimation by employing a result of Koukoulopoulos on the size of sifted Dirichlet $L$-functions (see Lemma \ref{fdp} below). 

The process to bound $S_1$ is not as much of a challenge and there are several ways to do it. Let us describe one. Thanks to the non-negativity of $\lambda$, we first ignore the terms $\Lambda(\ell)$ and bound them all crudely by $\log x$. Then, we count $z$-rough integers in the innermost summation, and conclude with an application of Lemma \ref{joka} which embodies the M\"{o}bius-mimicking behavior of the exceptional character on large primes. 

For $S_2$, we deal with the condition $P^-(k\ell)>z$ by using sieve weights again, and then, the new main idea is that we put Lemma \ref{compa} into action so that we switch from the sums of the convolution $1\ast\lambda$ over arithmetic progressions to its sums over integers that are coprime to $q$. In this procedure, we eventually reduce the estimation of $S_2$ to a bound that we already addressed in the estimation of $S_1$. 

All the above are but a brief explanation on how the right-hand of (\ref{tocon}) will be handled to reach a proof of Theorem \ref{maint}. The actual arguments and their relevant details are thoroughly developed in Sections \ref{1stpot} and \ref{lpot}.

\section{Preparatory Lemmas}

In this section, we state and prove auxiliary lemmas that are of use in the sequel of the text. They are grouped in four subsections based on their nature.

\subsection{Bounds on sums of $\lambda$ and $\lambda'$}
We start with a lemma on the sums of $\lambda=1\ast\chi$ over those integers that are coprime to a given $q\in\N$.

\begin{lemma}\label{csl}
Let $x\geq 1,\eps>0,q\in\N$, and $D$ be a positive integer such that $D\geq 2$. If $\chi$ is a non-principal Dirichlet character $\mod{D}$, then
\begin{align*}
\sum_{\substack{n\leq x\\(n,q)=1}}\lambda(n)=\frac{x\phi(q)}{q}\prod_{p\mid q}\bigg(1-\frac{\chi(p)}{p}\bigg)L(1,\chi)+O_{\eps}(\tau(q)x^{1/2}D^{1/2+\eps}).
\end{align*}
\end{lemma}
\begin{proof}
The proof is essentially a routine application of the hyperbola method, but for the sake of completeness, we provide a detailed derivation here.

Since $\sum_{m\leq u,(m,q)=1}1=u\phi(q)/q+O(\tau(q))$ for all $u\geq 1$ (this is a well-known consequence of M\"{o}bius inversion), we have that
\begin{align}\label{afhm}
\begin{split}
\sum_{\substack{n\leq x\\(n,q)=1}}\lambda(n)&=\sum_{\substack{k\leq \sqrt{x}\\(k,q)=1}}\chi(k)\sum_{\substack{m\leq x/k\\(m,q)=1}}1+\sum_{\substack{m\leq \sqrt{x}\\(m,q)=1}}\sum_{\substack{\sqrt{x}<k\leq x/m\\(k,q)=1}}\chi(k)\\
&=\frac{x\phi(q)}{q}\sum_{\substack{k\leq \sqrt{x}\\(k,q)=1}}\frac{\chi(k)}{k}+\sum_{\substack{m\leq \sqrt{x}\\(m,q)=1}}\sum_{\substack{\sqrt{x}<k\leq x/m\\(k,q)=1}}\chi(k)+O(\sqrt{x}).
\end{split}
\end{align}

Using M\"{o}bius inversion and the P\'{o}lya-Vinogradov inequality, we find that
\begin{align}\label{chM}
\sum_{\substack{k\leq u\\(k,q)=1}}\chi(k)=\sum_{k\leq u}\chi(k)\sum_{d\mid(k,q)}\mu(d)=\sum_{d\mid q}\mu(d)\chi(d)\sum_{\ell\leq u/d}\chi(\ell)\ll_{\eps}\tau(q)D^{1/2+\eps}\quad(u\geq 1).
\end{align}

Partial summation, combined with this bound, allow us to establish the convergence of the series $\sum_{k\geq 1,(k,q)=1}\chi(k)/k$ and assess the contribution of its tails. In particular,
\begin{align}\label{tb}
\sum_{\substack{k\leq \sqrt{x}\\(k,q)=1}}\frac{\chi(k)}{k}=\sum_{\substack{k\geq 1\\(k,q)=1}}\frac{\chi(k)}{k}+O_{\eps}\bigg(\frac{\tau(q)D^{1/2+\eps}}{\sqrt{x}}\bigg).
\end{align}

Inserting (\ref{tb}) into (\ref{afhm}), and applying the bound (\ref{chM}) twice to estimate the innermost character sum at the bottom line of (\ref{afhm}), we infer that
\begin{align*}
\sum_{\substack{n\leq x\\(n,q)=1}}\lambda(n)=\frac{x\phi(q)}{q}\sum_{\substack{k\geq 1\\(k,q)=1}}\frac{\chi(k)}{k}++O_{\eps}(\tau(q)x^{1/2}D^{1/2+\eps}).
\end{align*}

In order to conclude the proof, it suffices to show that the series $\sum_{k\geq 1,(k,q)=1}\chi(k)/k$ is equal to $L(1,\chi)\prod_{p\mid q}(1-\chi(p)/p)$. This is in fact relatively easy to see. Using the Euler product expansions of $L(\sigma,\chi)$ and $\sum_{k\geq 1,(k,q)=1}\chi(k)/k^{\sigma}$ for $\sigma>1$, we have that
$$\sum_{\substack{k\geq 1\\(k,q)=1}}\frac{\chi(k)}{k^{\sigma}}=\prod_{p\nmid q}\bigg(1-\frac{\chi(p)}{p^{\sigma}}\bigg)^{-1}=\prod_{p\mid q}\bigg(1-\frac{\chi(p)}{p^{\sigma}}\bigg)L(\sigma,\chi),$$
and now, the desired equality follows by taking the limit as $\sigma\to 1^+$ on both sides. This completes the proof of the lemma. 
\end{proof}

\begin{lemma}\label{compa}
Let $x\geq 1$ and $\eps\in(0,1/100)$. Let also $a$ and $q$ be two positive integers such that $q\leq x^{58/115-\eps}$ and $(a,q)=1$. If $D\geq 2$ is an integer and $\chi$ is a Dirichlet character $\mod{D}$, then
\begin{align*}
\sum_{\substack{k\ell\leq x\\k\ell\equiv a\mod{q}}}\lambda(k)=\frac{1}{\phi(q)}\sum_{\substack{k\ell\leq x\\(k\ell,q)=1}}\lambda(k)+O_{\eps}\bigg(\frac{Dx^{1-\eps/6}}{q}\bigg).
\end{align*}
\end{lemma}
\begin{proof}
This is basically stated as (5.6) with the error term given in (5.7) at \cite[p.~2048]{fip}. There, the authors have the additional condition $k>D^2$, but the proof that they outline on \cite[p.~2047-2048]{fip} as a proper adaptation of their work in \cite{fip2} applies to our case, too.
\end{proof}

Apart from the last two lemmas regarding the function $\lambda$, we will also need the following result which is concerned with the sums of $\lambda'=\chi\ast\log$ on arithmetic progressions.
\begin{lemma}\label{l's}
Let $a$ and $q$ be two coprime positive integers. Suppose also that $\chi$ is a real primitive Dirichlet character modulo $D\geq 2$. For $\eps\in(0,1/2)$ and $x>qD$, we have that
\begin{align*}
\sum_{\substack{n\leq x\\n\equiv a\mod{q}}}\lambda'(n)=\frac{1-\mathds{1}_{D\mid q}\chi(a)}{\phi(q)}\sum_{\substack{n\leq x\\(n,q)=1}}\lambda'(n)&+\frac{\mathds{1}_{D\mid q}\chi(a)}{q}(x\log x-x)\prod_{p\mid q}\bigg(1-\frac{\chi(p)}{p}\bigg)L(1,\chi)\\
&+O_{\eps}\bigg(\frac{x^{1/2+\eps}D^{1/2+\eps}}{q^{1/4}}\bigg).
\end{align*}
\end{lemma}
\begin{proof}
According to \cite[Proposition 4.2]{fip}, we have that
\begin{align}\label{ffi}
\begin{split}
&\sum_{\substack{n\leq x\\n\equiv a\mod{q}}}\lambda'(n)=\frac{1-\mathds{1}_{D\mid q}\chi(a)}{\phi(q)}\sum_{\substack{n\leq x\\(n,q)=1}}\lambda'(n)\\
&+\frac{\mathds{1}_{D\mid q}\chi(a)}{\phi(q)}\sum_{\substack{n\leq x\\(n,q)=1}}\lambda(n)\log n+O_{\eps}\bigg(\frac{D^{1/4}x^{1/2+\eps}}{q^{1/4}}\bigg).    
\end{split}
\end{align}

Moreover, combining Lemma \ref{csl} with partial summation, we obtain that
$$\sum_{\substack{n\leq x\\(n,q)=1}}\lambda(n)\log n=\frac{\phi(q)}{q}(x\log x-x)\prod_{p\mid q}\bigg(1-\frac{\chi(p)}{p}\bigg)L(1,\chi)+O_{\eps}(\tau(q)x^{1/2+\eps}D^{1/2+\eps}).$$

Plugging this asymptotic formula into (\ref{ffi}), and noticing that $\phi(q)/\tau(q)\!\gg\!q/(\tau(q)\log\log q)\!\gg_{\eps}\!q^{1-\eps}\geq q^{1/4}$, we conclude the proof of the current lemma.
\end{proof}

\subsection{Bounds related to sums of additive and multiplicative characters}
The next two lemmas constitute vital inputs in the forthcoming study of the sums $\sum_{n\leq x,n\equiv a\mod{q}}\lambda(n)$ and the proof of Proposition \ref{slaps} in Section \ref{sl'sec}.
\begin{lemma}\label{csaps}
Let $\chi$ be a primitive Dirichlet character $\mod{D}$. Assume that $b,M$ and $N$ are three positive integers and let $D'\in\N$ be a proper divisor of $D$. Then
\begin{align*}
\sum_{M\leq m<M+N}\chi(mD'+b)\ll \frac{D}{D'}.   
\end{align*}
\end{lemma}
\begin{proof}
Set $D^*=D/D'$. Since $D'$ is a proper divisor of $D$, it follows that $D^*>1$. Due to the $D$-periodicity of the character and the fact that its values are bounded by $1$, the proof reduces to showing that
\begin{align}\label{tpcl}
\sum_{j=1}^{D^*}\chi(jD'+b)=0.
\end{align}

Since $\chi$ is primitive, it is known that $\chi(n)=G(\bar{\chi})^{-1}\sum_{r=1}^D\bar{\chi}(r)e(rn/D)$, where $G(\bar{\chi})$ is the Gauss sum of the conjugate character $\bar{\chi}$. Hence,
\begin{align}\label{dtpc}
\sum_{j=1}^{D^*}\chi(jD'+b)=\frac{1}{G(\bar{\chi})}\sum_{r=1}^D\bar{\chi}(r)e\bigg(\frac{rb}{D}\bigg)\sum_{j=1}^{D^*}e\pthbgg{\frac{rj}{D^*}}.
\end{align}
But,
\begin{equation*}
\sum_{j=1}^{D^*}e\pthbgg{\frac{rj}{D^*}}=
\begin{cases}
D^*& \text{if }D^*\mid r,\\
0& \text{if }D^*\nmid r,
\end{cases}
\end{equation*}
and so (\ref{dtpc}) implies that
\begin{align}\label{tfcl}
\sum_{j=1}^{D^*}\chi(jD'+b)=\frac{D^*}{G(\bar{\chi})}\sum_{\substack{r=1\\r\equiv 0\mod{D^*}}}^D\bar{\chi}(r)e\pthbgg{\frac{rb}{D}}=\frac{D^*}{G(\bar{\chi})}\sum_{\ell=1}^{D'}\bar{\chi}(\ell D^*)e\pthbgg{\frac{\ell b}{D'}}.
\end{align}

Now, observe that $\chi(\ell D^*)=0$ for all $\ell\in\{1,\ldots,D'\}$ because $(\ell D^*,D)\geq D^*>1$. Employing this observation in (\ref{tfcl}), we deduce (\ref{tpcl}) and complete the proof of the lemma.
\end{proof}

We will also need the Weil bound on Kloosterman sums. Given natural numbers $m,n,c,$ the Kloosterman sum $K(m,n;c)$ is defined as
\begin{equation*}
K(m,n;c)\vcentcolon=\sum_{\substack{r=1\\(r,c)=1}}^ce\pthbgg{\frac{mr+n\bar{r}}{c}},
\end{equation*}
where $\bar{r}$ denotes the inverse of $r\mod{c}$.
\begin{lemma}[\textbf{The Weil bound on Kloosterman sums}]\label{wb}
Let $m,n$ and $c$ be positive integers. We have that
$$\abs{K(m,n;c)}\leq \tau(c)c^{1/2}(m,n,c)^{1/2},$$
where $(m,n,c)$ denotes the greatest common divisor of $m,n$ and $c$.
\end{lemma}

\subsection{Sieve lemmas}
We continue the current section with the statement and the proof of a sifting lemma. Its first part provides an explicit approximation of $\mathds{1}_{P^-(\cdot)>z}$ by a Dirichlet convolution $1\ast w$, where $w$ is a beta sieve weight. The first part of the lemma allows us to ``comfortably" interchange between $\mathds{1}_{P^-(\cdot)>z}$ and an appropriate sieve-weight convolution $1\ast w$ by introducing an error that can be sufficiently handled only with upper bounds. The second part of the lemma is a variant of the fundamental lemma of sieve theory containing logarithmic weights. The first part is essentially \cite[Lemma 3.2(i)]{mm}. The proof of the second part utilizes a result from the proof of \cite[Claim 3.9, p.~2946]{mine} and presents similarities with the proof of \cite[Lemma 3.2(ii)]{mm}.
\begin{lemma}\label{FL}
Let $z>1$. Furthermore, given an arbitrary $y\geq 1$, let $P(y)\vcentcolon=\prod_{p\leq y}p$. There exists $u_0\geq 4$ such that for every $u\geq u_0$, all of the following hold.
\begin{enumerate}[(a)]
\item For $r\in\N$, define
\begin{equation}\label{zr}
z_r\vcentcolon=z^{((u-2)/u)^r}.
\end{equation}
There exists an arithmetic function $w$ such that:
\vspace{1mm}
\begin{itemize}
\item $\abs{w(n)}\leq 1\,\text{ for all }\,n\in\N;$
\vspace{1.75mm}
\item $\supp(w)\subseteq \{d\in\N:d\mid \prod_{p\leq z}p\text{ and } d\leq z^u\},\text{ and};$
\item $\displaystyle{\mathds{1}_{P^-(n)>z}=(1\ast w)(n)+O\Big(\tau(n)^2\sum_{r\geq u/2}\mathds{1}_{P^-(n)>z_r}2^{-r}\Big)}\,
\text{ for all }\,n\in\N$.
\end{itemize}
\item Let $B\geq 1$ and $w$ be the arithmetic function of part (a). If $j$ is a non-negative integer and $\nu$ is a multiplicative function such that $\abs{\nu(p)}<\min\{B,p\}$ for every prime $p$, then if $u\geq 50B$, we have
\begin{align*}
\sum_{d\mid P(z)}\frac{w(d)\nu(d)(\log d)^j}{d}\ll_j(\log z)^j\prod_{p\leq z}\bigg(1-\frac{\nu(p)}{p}\bigg).
\end{align*}
\end{enumerate}
\end{lemma}
\begin{proof}
Part (a) follows from \cite[Lemma 3.2(i)]{mm} by letting our $u$ here play the role of $u\theta$ for a fixed $\theta\in(0,1/3)$, say $\theta=1/4$, and by taking $A=1,\beta_1=u_0/2$, and $\beta=u/2$.

For part (b), we give a detailed proof. To start, we first describe how the sieve weight $w$ of part (a) is constructed. Letting 
$$\CD\vcentcolon=\{1\}\cup\{d=p_1\cdots p_k\mid P(z):p_1>\cdots>p_k,\,p_1\cdots p_hp_h^{u/2}<z^u\text{ for all odd }h\leq k\},$$
the arithmetic function $w$ is defined as $w(d)=\mu(d)\mathds{1}_{\CD}(d)$ for $d\in\N$. In other words, we define it as an upper bound beta sieve weight.

Now, note that if $d=p_1\cdots p_k\mid P(z)$, and $p_1>\cdots>p_k$, but $d\notin \CD$, the last condition describing $\CD$ is violated at some positive odd integer $r\leq k$ for the first time, which means that there exists an odd integer $r$ such that $p_1\cdots p_rp_r^{u/2}\geq z^u$, while $p_1\cdots p_hp_h^{u/2}<z^{u}$ for all positive odd integers $h<r$. According to \cite[Section 6.3]{iwko}, for such a positive odd integer $r$, we have that $p_r\geq z_r$, where $z_r$ is defined in (\ref{zr}). Using the definition of $w$ along with all these, we obtain that
\begin{align}\label{bdef}
\begin{split}
&\absBgg{\sum_{d\mid P(z)}\frac{w(d)\nu(d)(\log d)^j}{d}-\sum_{d\mid P(z)}\frac{\mu(d)\nu(d)(\log d)^j}{d}}=\absBgg{\sum_{\substack{d\mid P(z)\\d\notin \CD}}\frac{\mu(d)\nu(d)(\log d)^j}{d}}\\
\leq\,&\sum_{r\text{ odd}}\absBgg{\sum_{\substack{p_r<\cdots<p_1\leq z,\\p_1\cdots p_rp_r^{u/2}\geq z^u,\\p_1\cdots p_hp_h^{u/2}<z^u\\\text{for all odd }h<r}}\frac{\mu(p_1\cdots p_r)\nu(p_1\cdots p_r)}{p_1\cdots p_r}\sum_{d'\mid P(p_r-1)}\frac{\mu(d')\nu(d')\log(d'p_1\cdots p_r)}{d'}}\\
\leq\,&\sum_{r\geq 1}\sum_{z_r\leq p_r<\cdots<p_1\leq z}\frac{\abs{\nu(p_1\cdots p_r)}}{p_1\cdots p_r}\absBgg{\sum_{d'\mid P(p_r-1)}\frac{\mu(d')\nu(d')\log(d'p_1\cdots p_r)}{d'}}.
\end{split}
\end{align}

Arguing as in the proof of \cite[(3.3), p.~2947]{mine}, one can show that
\begin{align}\label{1''}
\begin{split}
\sum_{d'\mid P(p_r-1)}\frac{\mu(d')\nu(d')\log(d'p_1\cdots p_r)}{d'}&\ll_j (\log z)^j(p_1\cdots p_r)^{\frac{1}{\log z}}\prod_{p<p_r}\bigg(1-\frac{\nu(p)}{p}\bigg)\\
&\leq (\log z)^je^r\prod_{p<z_r}\bigg(1-\frac{\nu(p)}{p}\bigg),
\end{split}
\end{align}
where the last line is justified by the inequalities $z_r\leq p_r\leq z$ and the assumption $\abs{\nu(p)}<p$. Similarly,
\begin{align}\label{2''}
\sum_{d\mid P(z)}\frac{\mu(d)\nu(d)(\log d)^j}{d}\ll_j (\log z)^j\prod_{p\leq z}\bigg(1-\frac{\nu(p)}{p}\bigg).
\end{align}

Let us now set $V(z)\vcentcolon=\prod_{p\leq z}(1-\nu(p)/p)$ and $S\vcentcolon=\sum_{d\mid P(z)}w(d)\nu(d)(\log d)^j/d$ for notational simplicity. We plug the bounds (\ref{1''}) and (\ref{2''}) into (\ref{bdef}) and infer that
\begin{align}\label{S}
\begin{split}
S&\ll_j (\log z)^jV(z)\Bigg(1+\sum_{r\geq 1}e^r\prod_{z_r\leq p\leq z}\bigg(1-\frac{\nu(p)}{p}\bigg)^{-1}\sum_{\substack{P^+(m)\leq z\\P^-(m)\geq z_r\\\omega(m)=r}}\frac{\mu(m)^2\abs{\nu(m)}}{m}\Bigg)\\
&\leq(\log z)^jV(z)\Bigg(1+\sum_{r\geq 1}2^{-r}\prod_{z_r\leq p\leq z}\bigg(1-\frac{\nu(p)}{p}\bigg)^{-1}\sum_{\substack{P^+(m)\leq z\\P^-(m)\geq z_r}}\frac{\mu(m)^2\abs{\nu(m)}6^{\omega(m)}}{m}\Bigg)\\
&\leq(\log z)^jV(z)\Bigg(1+\sum_{r\geq 1}2^{-r}\prod_{z_r\leq p\leq z}\bigg(1-\frac{\nu(p)}{p}\bigg)^{-1}\bigg(1+\frac{6\abs{\nu(p)}}{p}\bigg)\Bigg).
\end{split}
\end{align}

Since $\abs{\nu(p)}<\min\{B,p\}$, we have that $1+6\abs{\nu(p)}/p<e^{6B/p}$ and that $-\log(1-\nu(p)/p)<B/p+O_B(p^{-2})$. Combining these inequalities with Mertens' estimate, and recalling the definition of the $z_r$'s from (\ref{zr}), we conclude that
\begin{align*}
\prod_{z_r\leq p\leq z}\bigg(1-\frac{\nu(p)}{p}\bigg)^{-1}\bigg(1+\frac{6\abs{\nu(p)}}{p}\bigg)\ll \bigg(\frac{\log z}{\log z_r}\bigg)^{\!7B}\!=\bigg(\frac{u}{u-2}\bigg)^{7Br}.
\end{align*}
Inserting this bound into (\ref{S}), we deduce that
\begin{align*}
S\ll_j(\log z)^jV(z)\sum_{r\geq 1}\bigg(\frac{u^{7B}}{2(u-2)^{7B}}\bigg)^r.
\end{align*}

To complete the proof, it suffices to show that the series of the right-hand side is convergent. To this end, we have to show that for $u\geq 50B$, the fraction $u^{7B}/(2(u-2)^{7B})$ is smaller than a constant $c<1$. The proof will then be complete. Since $e^t-1\geq t$ for all $t\in\R$ and $\log 2>2/3$, it follows that $2^{1/(14B)}-1>(\log 2)/(14B)>1/(21B)$. Hence,
\begin{align*}
u\geq 50B>21B\cdot 2^{15/14}\geq 21B\cdot2^{\frac{14B+1}{14B}}>\frac{2^{\frac{14B+1}{14B}}}{2^{1/(14B)}-1}.
\end{align*}
But, the inequality $u>2^{(14B+1)/(14B)}/(2^{1/(14B)}-1)$ that we arrived at is equivalent to the inequality $u^{7B}/(2(u-2)^{7B})<2^{-1/2}$. Therefore, due to what was explained above, the proof of part (b), as well as that of the lemma, are complete.
\end{proof}

We can now use Lemma \ref{FL}(a) to prove the following sifting result.
\begin{lemma}\label{auxsf}
Let $f_1,f_2$ and $f_3$ be three given arithmetic functions such that $\abs{f_1(n)}\leq \log n$ and $\abs{f_2(n)},\abs{f_3(n)}\leq 1$ for all $n\in\N$. Let $x\geq z>1$ and set $f\vcentcolon=f_1\ast f_2\ast f_3$. Assume also that $a$ and $q$ are positive coprime integers such that $q\leq x^{1-\eps}$ for some $\eps\in(0,1)$. If $u_0$ is the real number from the statement of Lemma \ref{FL}, then for $u\geq \max\{50,u_0\}$, we have
\begin{align*}
\sum_{\substack{n\leq x\\n\equiv a\mod{q}\\P^-(n)>z}}f(n)=&\sum_{\substack{d_1,d_2,d_3\leq z^u\\(d_1d_2d_3,q)=1}}w(d_1)w(d_2)w(d_3)\sum_{\substack{k\ell m\leq \frac{x}{d_1d_2d_3}\\k\ell m\equiv a\bar{d_1d_2d_3}\mod{q}}}f_1(kd_1)f_2(\ell d_2)f_3(md_3)\\
&+O_{\eps}\bigg(e^{-u/15}\frac{x}{\phi(q)}\bigg(\frac{\log x}{\log z}\bigg)^{\!\!12}\,\bigg),
\end{align*}
where $\bar{d_1d_2d_3}$ denotes the inverse of $d_1d_2d_3\mod{q}$ and $w$ is the arithmetic function from Lemma \ref{FL}. If we have the stronger assumption $\abs{f_1(n)}\leq 1$ for all $n\in\N$, then the same holds with the only difference being a logarithmic power of $(\log x)^{11}$ in the big-Oh term instead of $(\log x)^{12}$.
\end{lemma}
\begin{proof}
We prove the lemma in the case where $\abs{f_1(n)}\leq \log n$ for all $n\in\N$. The proof of the other case is similar.

We first observe that
\begin{align}\label{twtp}
\begin{split}
&\sum_{\substack{d_1,d_2,d_3\leq z^u\\(d_1d_2d_3,q)=1}}w(d_1)w(d_2)w(d_3)\sum_{\substack{k\ell m\leq \frac{x}{d_1d_2d_3}\\k\ell m\equiv a\bar{d_1d_2d_3}\mod{q}}}f_1(kd_1)f_2(\ell d_2)f_3(md_3)\\
=&\sum_{\substack{k'\ell' m'\leq x\\k'\ell' m'\equiv a\mod{q}}}(1\ast w)(k')f_1(k')(1\ast w)(\ell')f_2(\ell')(1\ast w)(m')f_3(m').
\end{split}
\end{align}

We are going to make use of Lemma \ref{FL}(a) in order to replace the values of $1\ast w$ at the bottom line with the respective values of $\mathds{1}_{P^-(\cdot)>z}$. From the definition of the $z_r$'s, given in (\ref{zr}), we notice that $z>z_r$ for all $r\in\N$. Furthermore, we have that $\sum_{r\geq u/2}2^{-r}\asymp 2^{-u/2}$. Consequently, it follows that $\mathds{1}_{P^-(n)>z}\ll 2^{u/2}\sum_{r\geq u/2}2^{-r}\mathds{1}_{P^-(n)>z_r}$ for all $n\in\N$. This means that if we approximate each copy of $1\ast w$ in (\ref{twtp}) by $\mathds{1}_{P^-(\cdot)>z}$ according to Lemma \ref{FL}(a), then each occurrence of $\mathds{1}_{P^-(\cdot)>z}$ in the big-Oh term may be replaced with $2^{u/2}\sum_{r\geq u/2}2^{-r}\mathds{1}_{P^-(\cdot)>z_r}$. This way, we get
\begin{align}\label{a3a}
\begin{split}
&\sum_{\substack{k'\ell' m'\leq x\\ k'\ell' m'\equiv a\mod{q}}}(1\ast w)(k')f_1(k')(1\ast w)(\ell')f_2(\ell')(1\ast w)(m')f_3(m')\\
=&\sum_{\substack{n\leq x\\n\equiv a\mod{q}\\P^-(n)>z}}f(n)+O\Big(2^u(\log x)\sum_{r_1,r_2,r_3\geq u/2}2^{-(r_1+r_2+r_3)}\sum_{\substack{n\leq x\\n\equiv a\mod{q}}}g_{r_1,r_2,r_3}(n)\Big),
\end{split}
\end{align}
where $g_{r_1,r_2,r_3}$ are the multiplicative functions defined as the following triple convolutions:
$$g_{r_1,r_2,r_3}(n)\vcentcolon=\sum_{k\ell m=n}\mathds{1}_{P^-(k)>z_{r_1}}\tau(k)^2\mathds{1}_{P^-(\ell)>z_{r_2}}\tau(\ell)^2\mathds{1}_{P^-(m)>z_{r_3}}\tau(m)^2\quad (n\in\N).$$

We are now going to bound the sum $\sum_{n\leq x,n\equiv a\mod{q}}g_{r_1,r_2,r_3}(n)$ when $z_{r_1}\leq z_{r_2}\leq z_{r_3}$. The bound of the sum in all the other possible orderings, such as when $z_{r_2}\leq z_{r_3}\leq z_{r_1}$, is the same and can be proven similarly. For this reason, we only cover in detail the case where $z_{r_1}\leq z_{r_2}\leq z_{r_3}$.

If $d\mid n,$ then obviously $\mathds{1}_{P^-(d)>z_{r_i}}\tau(d)^2\leq \tau(n)^2$ for $i\in\{1,2,3\}$. Therefore, the non-negative multiplicative function $g_{r_1,r_2,r_3}$ satisfies the inequality $g_{r_1,r_2,r_3}(n)\leq \tau(n)^6\tau_3(n)$ for all $n\in\N$. This inequality allows us to verify that the celebrated theorem of Shiu \cite[Theorem 1]{sh} is applicable to $g_{r_1,r_2,r_3}$. Making use of it in conjunction with Mertens' theorem and the definition (\ref{zr}), we get
\begin{align*}
\sum_{\substack{n\leq x\\n\equiv a\mod{q}}}g_{r_1,r_2,r_3}(n)&\ll_{\eps}\frac{x}{\phi(q)\log x}\exp\bigg(\sum_{z_{r_1}<p\leq z_{r_2}}\frac{4}{p}+\sum_{z_{r_2}<p\leq z_{r_3}}\frac{8}{p}+\sum_{z_{r_3}<p\leq x}\frac{12}{p}\bigg)\\
&\ll \frac{x}{\phi(q)\log x}\bigg(\frac{\log z_{r_2}}{\log z_{r_1}}\bigg)^4\bigg(\frac{\log z_{r_3}}{\log z_{r_2}}\bigg)^8\bigg(\frac{\log x}{\log z_{r_3}}\bigg)^{12}\\
&=\frac{x(\log x)^{11}}{\phi(q)(\log z_{r_1})^4(\log z_{r_2})^4(\log z_{r_3})^4}=\frac{x(\log x)^{11}}{\phi(q)(\log z)^{12}}\bigg(\frac{u}{u-2}\bigg)^{4(r_1+r_2+r_3)},
\end{align*}
where at the last step above, we recalled the definition of the $z_r$'s from (\ref{zr}).

As we mentioned, the bound that we just obtained can be similarly established for the sums of $g_{r_1,r_2,r_3}$ when one of the other possible orderings of $z_{r_1},z_{r_2}$ and $z_{r_3}$ holds. So, from (\ref{a3a}), it follows that
\begin{align}\label{amsubs}
\begin{split}
&\sum_{\substack{k'\ell' m'\leq x\\k'\ell' m'\equiv a\mod{q}}}(1\ast w)(k')f_1(k')(1\ast w)(\ell')f_2(\ell')(1\ast w)(m')f_3(m')\\
=&\sum_{\substack{n\leq x\\n\equiv a\mod{q}\\P^-(n)>z}}f(n)+O_{\eps}\Bigg(\frac{2^ux(\log x)^{12}}{\phi(q)(\log z)^{12}}\prod_{i=1}^3\sum_{r_i\geq u/2}\bigg(\frac{u^4}{2(u-2)^4}\bigg)^{\!r_i}\bigg).
\end{split}
\end{align}

Now, since $u\geq 50$, it turns out that $u^4/(2(u-2)^4)\leq2^{-3/4}$. Then
$$\sum_{r_i\geq u/2}\bigg(\frac{u^4}{2(u-2)^4}\bigg)^{\!r_i}\leq \sum_{r_i\geq u/2}2^{-3r_i/4}\ll 2^{-3u/8}\quad(i=1,2,3),$$
and so (\ref{amsubs}) yields that
\begin{align}\label{ending}
\begin{split}
&\sum_{\substack{k'\ell' m'\leq x\\k'\ell' m'\equiv a \mod{q}}}(1\ast w)(k')f_1(k')(1\ast w)(\ell')f_2(\ell')(1\ast w)(m')f_3(m')\\
=&\sum_{\substack{n\leq x\\n\equiv a\mod{q}\\P^-(n)>z}}f(n)+O_{\eps}\bigg(2^{-u/8}\frac{x}{\phi(q)}\bigg(\frac{\log x}{\log z}\bigg)^{\!\!12}\,\bigg).
\end{split}
\end{align}

Since $2^{-u/8}<e^{-u/15}$, we complete the proof by substituting the bottom line of (\ref{twtp}) into the bottom line of (\ref{ending}).
\end{proof}

Upon a vigilant examination of the proof of Lemma \ref{auxsf}, we observe that the lemma can be generalized for convolutions of arithmetic functions whose values are bounded by expressions of the form $(\log n+1)^{A}\tau(n)^B$ for some $A,B>0$. We decided to state and prove Lemma \ref{auxsf} for triple convolutions under weaker size assumptions on $f_1,f_2$ and $f_3$ because its given formulation is precisely what we need later. 

Lemma \ref{auxsf} proves to be useful when we wish to study sifted sums $\sum_{n\leq x,n\equiv a\mod{q},P^-(n)>z}f(n)$, where $f$ is not completely multiplicative but can be written as a triple convolution of completely multiplicative or completely additive functions. Indeed, if we employed a sieve-weight convolution $1\ast w$ to approximate the original sum by $\sum_{(d,q)=1}w(d)\sum_{m\leq x/d,m\equiv a\bar{d}\mod{q}}f(md)$, then we could not ``untangle" the variables $m$ and $d$. This is an ``inconvenient" technicality. However, if we resort to Lemma \ref{auxsf}, this issue is circumvented because the variables can be separated in the arguments of the ``components" $f_1$, $f_2$ and $f_3$ which are either completely multiplicative or completely additive.

\subsection{Results on the M\"{o}bius-like behavior of exceptional characters} 
The next result is a rigorous way of phrasing that the character of a Landau-Siegel zero of large quality tends to ``behave" like the M\"{o}bius function on large primes.
\begin{lemma}\label{joka}
Let $\chi$ be a primitive quadratic character modulo $D\geq 2$. Assume also that $L(\cdot,\chi)$ has a real zero $\beta$ such that $\beta=1-1/(\eta\log D)$ for some $\eta\geq 10$. If $x=D^V$ for some $V>2$, and
$$z\vcentcolon=D^{\min\{\sqrt{V/(\log \eta)},\,2\}},$$
then there exists some absolute constant $c>0$ such that
\begin{eqnarray*}
\sum_{\substack{z<n\leq x\\P^-(n)>z}}\frac{(1\ast\chi)(n)}{n}\ll\frac{V^3}{\eta}+\exp\Big(\!\!-c\sqrt{V\log \eta}\Big).
\end{eqnarray*}    
\end{lemma}
\begin{proof}
We first set $v=\min\{\sqrt{V/(\log \eta)},2\}$. Then, according to \cite[Lemma 2.2]{mm}, we have that
\begin{equation}\label{fjkp}
\sum_{\substack{z<n\leq x\\P^-(n)>z}}\frac{(1\ast\chi)(n)}{n}\ll\frac{V^3}{v^3}\pthbgg{\frac{1}{v^2\eta^{v/2}}+\frac{v}{\eta}+\frac{1}{z}}.
\end{equation}

If $v=2$, then (\ref{fjkp}) yields the bound $V^3/\eta+V^3/D^2$. However, Siegel's estimate (\ref{siegel}) implies that $\eta\ll D^2$, and so in the case $v=2$, the bound is $V^3/\eta$. This concludes the proof of the lemma when $v=2$. 

If $V\in(2,4\log\eta)$, or equivalently $v=\sqrt{V/(\log\eta)}<2$, then the middle term in the parentheses of the right-hand side of (\ref{fjkp}) is absorbed by the first term. Consequently, using (\ref{fjkp}) along with the bound $\log\eta\ll\log D$, which follows from (\ref{siegel}), we infer that there exists an absolute constant $C\in(0,1]$ such that
\begin{align*}
\sum_{\substack{z<n\leq x\\P^-(n)>z}}\frac{(1\ast\chi)(n)}{n}&\ll (V\log \eta)^{5/2}\bigg(\exp\Big(\!\!-\sqrt{V\log\eta}/2\Big)+\exp\bigg(\!\!-\sqrt{\frac{V(\log D)^2}{\log\eta}}\bigg)\bigg)\\
&\ll (V\log\eta)^{5/2}\exp\Big(\!\!-C\sqrt{V\log\eta}\Big)\ll \exp\Big(\!\!-c\sqrt{V\log\eta}\Big),
\end{align*}
where $c\vcentcolon=C/2$. We have thus completed the proof of the current lemma in both cases for $v$.
\end{proof}

We close the section with a lemma which will be necessary for the estimation of the difference $\Delta$ defined in (\ref{Delta}). It can be extracted from \cite[Exercise 22.3(e), p.~232]{dimb}, but since the solution to that exercise is not included in the cited book, we provide a proof which is based off a theorem that is proven in the work of Koukoulopoulos on multiplicative functions with small averages in \cite{oldk}. 
\begin{lemma}\label{fdp}
Let $D\geq 2$ be an integer and $\chi$ be a real primitive Dirichlet character $\mod{D}$. If $L(\cdot,\chi)$ has a Landau-Siegel zero $\beta=1-1/(\eta\log D)$ for some $\eta\geq 2$ and $z\in(1,D^2]$, then
$$L(1,\chi)\prod_{p\leq z}\bigg(1-\frac{\chi(p)}{p}\bigg)\ll \frac{\log D}{\eta\log z}.$$
\end{lemma}
\begin{proof}
In light of the P\'{o}lya-Vinogradov inequality, there exists some sufficiently large absolute constant $C\geq 2$ such that the main condition of \cite[Theorem 2.4]{oldk} (that is, assumption (1.10) on \cite[p.~1574]{oldk}) is met with $f=\chi$. Hence, according to \cite[Theorem 2.4]{oldk} (and following the notation of the cited paper), we deduce that
\begin{align}\label{2.4dp}
L(1,\chi)\prod_{p\leq D^C}\bigg(1-\frac{\chi(p)}{p}\bigg)\ll_C(1-\beta)\log D=\frac{1}{\eta}.
\end{align}
Since Mertens' theorem implies that
\begin{align*}
\prod_{p\leq D^C}\bigg(1-\frac{\chi(p)}{p}\bigg)&=\prod_{p\leq z}\bigg(1-\frac{\chi(p)}{p}\bigg)\prod_{z<p\leq D^C}\bigg(1-\frac{\chi(p)}{p}\bigg)\\
&\gg \prod_{p\leq z}\bigg(1-\frac{\chi(p)}{p}\bigg)\exp\bigg(\!\!-\sum_{z<p\leq D^C}\frac{\chi(p)}{p}\bigg)\\
&\geq \prod_{p\leq z}\bigg(1-\frac{\chi(p)}{p}\bigg)\exp\bigg(\!\!-\sum_{z<p\leq D^C}\frac{1}{p}\bigg)\\
&\gg_C\prod_{p\leq z}\bigg(1-\frac{\chi(p)}{p}\bigg)\frac{\log z}{\log D},
\end{align*}
the estimate (\ref{2.4dp}) gives $L(1,\chi)(\log z)(\log D)^{-1}\prod_{p\leq z}(1-\chi(p)/p)\ll \eta^{-1}$, which in turn yields the bound of the lemma by simply taking $\log z$ and $\log D$ to the other side.
\end{proof}

\section{The sums of $\lambda$ over arithmetic progressions}\label{sl'sec}

This section is exclusively dedicated to the proof of Proposition \ref{slaps} below.

Proposition \ref{slaps} furnishes an asymptotic formula for the sums of $\lambda$ on arithmetic progressions with an error that is better than $\tau(q)\sqrt{x}D^{1/2+\eps}$, which is what a rather typical hyperbola method approach would produce. This feature of Proposition \ref{slaps} is the reason of its inclusion in this paper since an error larger than $\sqrt{x}$ would not allow us to go beyond the square-root barrier of GRH. 

We decided to prove Proposition \ref{slaps} here separately because a proper treatment of the sums $\sum_{n\leq x,n\equiv a \mod{q}}\lambda(n)$ needs considerably more effort than the proofs of the previous section.

The proof of Proposition \ref{slaps} does not follow the Vorono\u{\i}-type routes. 
Instead, we adopt a ``\`{a} la Heath-Brown" approach in the sense that our method constitutes an adaptation of Heath-Brown's arguments from his work on the ternary divisor function in \cite{hbd3}. Upon comparing the quoted work of Heath-Brown to our proof below, one can realize that the case of $\tau_3$ presents more difficulties. This should come as no surprise given that we are dealing with the function $\lambda=1\ast\chi$ that basically mimics the classical divisor function $\tau$ which is of lower complexity than $\tau_3$.

\begin{proposition}\label{slaps}
Let $x>1,\eps\in(0,1/15)$, and suppose that $a$ and $q$ are coprime positive integers such that $q\leq x^{2/3-\eps}$. Let also $\chi$ be a real primitive Dirichlet character modulo $D\geq 2$. Then
\begin{align*}
\sum_{\substack{n\leq x\\n\equiv a\mod{q}}}\lambda(n)=\frac{x(1+\mathds{1}_{D\mid q}\chi(a))}{q}\prod_{p\mid q}\bigg(1-\frac{\chi(p)}{p}\bigg)L(1,\chi)+O_{\eps}\bigg(\frac{Dx^{1-\eps/2}}{q}\bigg).
\end{align*}
\end{proposition}

Since the proof of Proposition \ref{slaps} is long, we divide it in subsections. In the first subsection, we open the convolution $\lambda=1\ast\chi$, stabilize the value of the character over a given residue class, and restrict the summation variables to short intervals. This is the starting step of the proof.

However, before diving into the first step of the proof of Proposition \ref{slaps}, let us first mention that throughout the following subsections of Section \ref{sl'sec}, the symbols $\bar{m}$ and $m^*$ denote the multiplicative inverses of the positive integer $m$ $\mod{q}$ and $\mod{(D,q)}$, respectively.

\subsection{The preliminary transformation} We have that
\begin{align}\label{beg}
\sum_{\substack{n\leq x\\n\equiv a\mod{q}}}\lambda(n)=\sum_{\substack{r\mod{D}\\(r,D)=1}}\chi(r)\sum_{\substack{k\ell\leq x\\k\ell\equiv a\mod{q}\\k\equiv r\mod{D}}}1.
\end{align}

Now, let $\delta\vcentcolon=q^{1+\eps}x^{-1}\leq 1$, and consider the double short sums
\begin{align}\label{dNKL}
N(K,L;r)\vcentcolon=\sum_{\substack{K<k\leq (1+\delta)K\\L<\ell\leq (1+\delta)L\\k\ell\equiv a\mod{q}\\k\equiv r\mod{D}}}1.
\end{align}
The reason for considering these sums is to ``decouple" the variables $k$ and $\ell$ at the right-hand side of (\ref{beg}). Using the sums $N(K,L;r)$, we rewrite the inner sum from the right-hand side of (\ref{beg}) as
\begin{align}\label{sNs}
\sum_{\substack{k\ell\leq x\\k\ell\equiv a\mod{q}\\k\equiv r\mod{D}}}1=\sum_{KL\leq x}N(K,L;r)+O\bigg(\sum_{\frac{x}{(1+\delta)^2}<KL\leq x}N(K,L;r)\bigg),
\end{align}
where here, as well as in any other subsequent occurrence, in the sums over $K$ and $L$, both parameters $K$ and $L$ run over the set of powers $1+\delta,(1+\delta)^2,(1+\delta)^3,\ldots$.

We would like to bound the big-Oh contribution. Dropping the condition $k\equiv r\mod{D}$ from the sums $N(K,L;r)$, we obtain that
\begin{align}\label{gfS}
\sum_{\frac{x}{(1+\delta)^2}<KL\leq x}N(K,L;r)\leq \sum_{\substack{\frac{x}{(1+\delta)^2}<n\leq x(1+\delta)^2\\n\equiv a\mod{q}}}\tau(n).
\end{align}
Now, our aim is to apply Shiu's theorem for sums of divisor functions \cite[Theorem 2]{sh} (in the cited theorem, Shiu denotes the classical divisor function by $d_2$) to the divisor sum on the right. In order to do that, we have to verify that the length of the sum is at least $q^{1+\eps}$. It really is because
\begin{align*}
x(1+\delta)^2-\frac{x}{(1+\delta)^2}>x((1+\delta)^2-1)>\delta x=q^{1+\eps}.
\end{align*}
Hence, Shiu's theorem is applicable, and so combining it with Mertens' theorem and the bound $\phi(q)\gg q/(\log\log q)\gg_{\eps}q^{1-\eps}$, we derive that
\begin{align*}
\sum_{\substack{\frac{x}{(1+\delta)^2}<n\leq x(1+\delta)^2\\n\equiv a\mod{q}}}\tau(n)&\ll_\eps\frac{1}{\phi(q)\log x}\bigg(x(1+\delta)^2-\frac{x}{(1+\delta)^2}\bigg)\exp\bigg\{\sum_{p\leq 4x}\frac{2}{p}\bigg\}\\
&\ll\frac{x\log x}{\phi(q)(1+\delta)^2}\Big((1+\delta)^4-1\Big)\leq\frac{x\log x}{\phi(q)}\Big(\delta^4+4\delta^3+6\delta^2+4\delta\Big)\\
&\ll \frac{\delta x\log x}{\phi(q)}=\frac{q^{1+\eps}\log x}{\phi(q)}\ll_{\eps} x^{3\eps}\leq\frac{x^{1-\eps}}{q}.
\end{align*}
We plug the estimate that we just obtained into (\ref{gfS}), and we then insert the resulting bound into (\ref{sNs}). This way, we infer that
\begin{align*}
\sum_{\substack{k\ell\leq x\\k\ell\equiv a\mod{q}\\k\equiv r\mod{D}}}1=\sum_{KL\leq x}N(K,L;r)+O_{\eps}\bigg(\frac{x^{1-\eps}}{q}\bigg).
\end{align*}

Consequently, (\ref{beg}) yields
\begin{align}\label{aftdd}
\sum_{\substack{n\leq x\\n\equiv a\mod{q}}}\lambda(n)=\sum_{\substack{r\mod{D}\\(r,D)=1}}\chi(r)\sum_{KL\leq x}N(K,L;r)+O_{\eps}\bigg(\frac{Dx^{1-\eps}}{q}\bigg),
\end{align}
and so the proof of Proposition \ref{slaps} reduces to the study of the sums $N(K,L;r)$.

\subsection{Detection of the $\mod{q}$ congruence conditions of $N(K,L;r)$ via additive characters}
Recalling the definition of $N(K,L;r)$ from (\ref{dNKL}) and using the elementary geometric sum
\begin{equation*}
\sum_{j=1}^{q}e\pthbgg{\frac{mj}{q}}=
\begin{cases}
q& \text{if }q\mid m,\\
0& \text{if }q\nmid m,
\end{cases}
\end{equation*}
we have that
\begin{align}\label{long}
\begin{split}
N(K,L;r)&=\sum_{\substack{\beta,\gamma \mod{q}\\\beta\gamma\equiv a\mod{q}}}\sum_{\substack{K<k\leq (1+\delta)K\\L<\ell\leq (1+\delta)L\\k\equiv \beta\mod{q}\\k\equiv r\mod{D}\\\ell\equiv \gamma\mod{q}}}1\\
&=\frac{1}{q^2}\sum_{\substack{\beta,\gamma \mod{q}\\\beta\gamma\equiv a\mod{q}}}\sum_{\substack{K<k\leq (1+\delta)K\\L<\ell\leq (1+\delta)L\\k\equiv r\mod{D}}}\sum_{s,t=1}^qe\bigg(\frac{s(\beta-k)}{q}\bigg)e\bigg(\frac{t(\gamma-\ell)}{q}\bigg)\\
&=\frac{1}{q^2}\sum_{s,t=1}^qK(s,at;q)\Bigg(\sum_{\substack{K<k\leq (1+\delta)K\\k\equiv r\mod{D}}}e\bigg(-\frac{sk}{q}\bigg)\Bigg)\Bigg(\sum_{L<\ell\leq (1+\delta)L}e\bigg(-\frac{t\ell}{q}\bigg)\Bigg)\\
&=\frac{1}{q^2}\sum_{s,t=1}^qK(s,at;q)e\bigg(\!-\frac{sr}{q}\bigg)F(s)G(t),
\end{split}
\end{align}
where $K(s,at;q)$ is the Kloosterman sum
\begin{align*}
K(s,at;q)=\sum_{\substack{\beta=1\\(\beta,q)=1}}^qe\bigg(\frac{s\beta+at\bar{\beta}}{q}\bigg)=\sum_{\substack{\beta,\gamma\mod{q}\\\beta\gamma\equiv a\mod{q}}}e\bigg(\frac{s\beta+t\gamma}{q}\bigg),
\end{align*}
and
\begin{align}\label{FG}
F(s)\vcentcolon=\sum_{\frac{K-r}{D}<k'\leq \frac{(1+\delta)K-r}{D}}e\bigg(\!-\frac{k'sD}{q}\bigg),\quad\quad G(t)\vcentcolon=\sum_{L<\ell\leq (1+\delta)L}e\bigg(\!-\frac{\ell t}{q}\bigg).
\end{align}

Let $q'\vcentcolon=q/(D,q)$. At the bottom line of (\ref{long}), we separate the terms corresponding to $t=q$ and to $s=mq'$ for $m\in\{1,\ldots,(D,q)\}$. Upon noticing that $G(q)=G(0)$ and $F(mq')=F(0)$ for all $m\in\{1,\ldots,(D,q)\}$, we then write
\begin{align}\label{rewNKL}
\begin{split}
N(K,L;r)=&\,\frac{G(0)}{q^2}\sum_{s=1}^qK(s,0;q)e\bigg(\!-\frac{sr}{q}\bigg)F(s)\\
&+\frac{F(0)}{q^2}\sum_{m=1}^{(D,q)}\sum_{t=1}^qK(mq',at;q)e\bigg(\!-\frac{mr}{(D,q)}\bigg)G(t)\\
&+\frac{1}{q^2}\sum_{\kappa=0}^{D-1}\sum_{\substack{\frac{\kappa q}{D}<s<\frac{(\kappa+1)q}{D}\\0<t<q}}K(s,at;q)e\bigg(-\frac{sr}{q}\bigg)F(s)G(t).
\end{split}
\end{align}

The main contribution in $N(K,L;r)$ comes from the first two lines in (\ref{rewNKL}). The multiple sum of its bottom line turns out to yield negligible contribution.

\subsection{Analysis of $N(K,L;r)$ -- Part I: Small contributions}
In the current subsection, we are bounding the contribution of
\begin{align*}
E\vcentcolon=\frac{1}{q^2}\sum_{\kappa=0}^{D-1}\sum_{\substack{\frac{\kappa q}{D}<s<\frac{(\kappa+1)q}{D}\\0<t<q}}K(s,at;q)e\bigg(-\frac{sr}{q}\bigg)F(s)G(t).
\end{align*}
Our goal toward the estimation of these triple sums is to bound the inner sums over $s$ and $t$. In this direction, we need the basic exponential sum estimate $\sum_{X<n\leq X+Y}e(-\alpha n)\ll \rVert\alpha\lVert^{-1}$ which holds for all $X,Y\geq 0$ and every non-integer $\alpha\in\R$. Here $\rVert\alpha\lVert$ denotes the distance of $\alpha$ from the closest integer. Consulting the definitions of $F(s)$ and $G(t)$ from (\ref{FG}), and combining the quoted exponential sum estimate with Lemma \ref{wb} and the trivial inequality $(s,t,q)^{1/2}\leq \sum_{d\mid(s,t,q)}d^{1/2}$, we deduce that
\begin{align}\label{fec}
\sum_{\substack{\frac{\kappa q}{D}<s<\frac{(\kappa+1)q}{D}\\0<t<q}}K(s,at;q)e\bigg(-\frac{sr}{q}\bigg)F(s)G(t)&\ll \tau(q)\sqrt{q}\sum_{\substack{\frac{\kappa q}{D}<s<\frac{(\kappa+1)q}{D}\\0<t<q}}\frac{(s,t,q)^{1/2}}{\lVert t/q\rVert\cdot\lVert sD/q\rVert}\nn
&\leq\tau(q)\sqrt{q}\sum_{d\mid q}d^{1/2}\!\!\sum_{\substack{\frac{\kappa q}{D}<s<\frac{(\kappa+1)q}{D}\\0<t<q\\d\mid s,\,d\mid t}}\frac{1}{\lVert t/q\rVert\cdot\lVert sD/q\rVert}\nn
&=\tau(q)\sqrt{q}\sum_{d\mid q}d^{1/2}\!\!\sum_{\substack{\frac{\kappa q}{dD}<s'<\frac{(\kappa+1)q}{dD}\\0<t'<q/d}}\frac{1}{\lVert dt'/q\rVert\cdot\lVert dDs'/q\rVert}.
\end{align}

Since $\lVert u\rVert=\abs{u}$ when $\abs{u}\leq 1/2$, and $\lVert\cdot\rVert$ is even and $1$-periodic, we have that
\begin{align}\label{t'sum}
\begin{split}
&\sum_{0<t'<q/d}\frac{1}{\lVert dt'/q\rVert}=\frac{q}{2d}\sum_{0<t'\leq q/(2d)}\frac{1}{t'}+\sum_{q/(2d)<t'<q/d}\frac{1}{\lVert dt'/q\rVert}\\
\ll&\,\frac{q\log q}{d}+\sum_{0<\mathfrak{t}<q/(2d)}\frac{1}{\lVert d\mathfrak{t}/q\rVert}\ll\frac{q\log q}{d}+\frac{q}{d}\sum_{0<\mathfrak{t}<q/(2d)}\frac{1}{\mathfrak{t}}\ll\frac{q\log q}{d}.
\end{split}
\end{align}
In a similar fashion, one can show that
\begin{align}\label{s'sum}
\sum_{\frac{\kappa q}{dD}<s'<\frac{(\kappa+1)q}{dD}}\frac{1}{\lVert dDs'/q\rVert}\ll \frac{q\log q}{dD}.
\end{align}

We substitute the sums of (\ref{t'sum}) and (\ref{s'sum}) with their respective bounds into (\ref{fec}) and obtain
\begin{align*}
\sum_{\substack{\frac{\kappa q}{D}<s<\frac{(\kappa+1)q}{D}\\0<t<q}}K(s,at;q)e\bigg(-\frac{sr}{q}\bigg)F(s)G(t)\ll \frac{\tau(q)q^{5/2}(\log q)^2}{D}\sum_{d\mid q}\frac{1}{d^{3/2}}\ll \frac{\tau(q)q^{5/2}(\log q)^2}{D}.
\end{align*}

Hence,
\begin{align}\label{E}
E\ll \tau(q)\sqrt{q}(\log q)^2\ll\frac{x^{1-\eps/2}}{q},
\end{align}
where the last step can be justified by the divisor bound and the assumption $q\leq x^{2/3-\eps}$.

\subsection{Analysis of $N(K,L;r)$ -- Part II: The main contributions}

The first main contribution in (\ref{rewNKL}) is coming from the term
\begin{align*}
M_1\vcentcolon=\frac{G(0)}{q^2}\sum_{s=1}^qK(s,0;q)e\bigg(\!-\frac{sr}{q}\bigg)F(s).
\end{align*}
Recalling the definitions of $F$ and $G$ from (\ref{FG}), we can calculate that
\begin{align}\label{M1}
\begin{split}
M_1&=\frac{G(0)}{q^2}\sum_{\substack{K<k\leq (1+\delta)K\\k\equiv r\mod{D}\\1\leq \beta\leq q\\(\beta,q)=1}}\sum_{s=1}^qe\bigg(\frac{s(\beta-k)}{q}\bigg)=\frac{G(0)}{q}\sum_{\substack{\beta=1\\(\beta,q)=1}}^q\sum_{\substack{K<k\leq (1+\delta)K\\k\equiv \beta\mod{q}\\k\equiv r\mod{D}}}1\\
&=\frac{G(0)}{q}\sum_{\substack{K<k\leq (1+\delta)K\\k\equiv r\mod{D}\\(k,q)=1}}1=\frac{1}{q}\sum_{\substack{K<k\leq (1+\delta)K\\L<\ell\leq (1+\delta)L\\k\equiv r\mod{D}\\(k,q)=1}}1.
\end{split}
\end{align}

Next, we look at
\begin{align*}
M_2\vcentcolon=\frac{F(0)}{q^2}\sum_{m=1}^{(D,q)}\sum_{t=1}^qK(mq',at;q)e\bigg(\!-\frac{mr}{(D,q)}\bigg)G(t),
\end{align*}
which can also yield major contributions to the asymptotic evaluation of $N(K,L;r)$. Since $q'=q/(D,q)$ and $K(mq',at;q)=\sum_{\gamma=1,(\gamma,q)=1}^qe((mq'\bar{\gamma}+at\gamma)/q)$, it follows that
\begin{align}\label{M2c}
\begin{split}
M_2&=\frac{F(0)}{q^2}\sum_{m=1}^{(D,q)}e\bigg(\!-\frac{mr}{(D,q)}\bigg)\sum_{\substack{L<\ell\leq (1+\delta)L}}\sum_{\substack{\gamma=1\\(\gamma,q)=1}}^qe\bigg(\frac{m\bar{\gamma}}{(D,q)}\bigg)\sum_{t=1}^qe\bigg(\frac{t(a\gamma-\ell)}{q}\bigg)\\
&=\frac{F(0)}{q}\sum_{m=1}^{(D,q)}e\bigg(\!-\frac{mr}{(D,q)}\bigg)\sum_{L<\ell\leq (1+\delta)L}\sum_{\substack{\gamma=1\\(\gamma,q)=1\\a\gamma\equiv\ell\mod{q}}}^qe\bigg(\frac{m\bar{\gamma}}{(D,q)}\bigg)\\
&=\frac{F(0)}{q}\sum_{L<\ell\leq (1+\delta)L}\sum_{\substack{\gamma=1\\(\gamma,q)=1\\a\gamma\equiv\ell\mod{q}}}^q\sum_{m=1}^{(D,q)}e\bigg(\frac{m(\bar{\gamma}-r)}{(D,q)}\bigg)\\
&=\frac{F(0)(D,q)}{q}\sum_{\substack{L<\ell\leq (1+\delta)L\\(\ell,q)=1}}\sum_{\substack{\gamma=1\\\gamma\equiv\bar{a}\ell\mod{q}\\\gamma\equiv r^*\mod{(D,q)}}}^q1.
\end{split}
\end{align}
Notice that $r^*$, that is, the inverse of $r\mod{(D,q)}$, does exist because $r$ is coprime to $D$, and this can be seen from the introduction of the parameter $r$ in (\ref{beg}).

Now, thanks to the Chinese Remainder Theorem, there is a unique solution $\gamma\mod{q}$ to the system of linear congruences $\gamma\equiv\bar{a}\ell\mod{q}$ and $\gamma\equiv r^*\mod{(D,q)}$ if and only if $\bar{a}\ell\equiv r^*\mod{(D,q)}$, or equivalently, if and only if $r\ell\equiv a\mod{(D,q)}$. Furthermore, from the definition of $F(s)$ in (\ref{FG}), we note that $F(0)=\sum_{K<k\leq (1+\delta)K,k\equiv r\mod{D}}1$. Therefore, continuing from the last line of (\ref{M2c}), we conclude that
\begin{align}\label{M2}
M_2=\frac{(D,q)}{q}\sum_{\substack{K<k\leq (1+\delta)K\\L<\ell\leq (1+\delta)L\\k\equiv r\mod{D}\\(\ell,q)=1\\k\ell\equiv a\mod{(D,q)}}}1.
\end{align}

\subsection{Returning back to the sums of $\lambda$ -- Endgame: Completion of the proof of Proposition \ref{slaps}}
Let us recall that $M_1,M_2$ and $E$ were defined as the top, middle and bottom line of (\ref{rewNKL}), respectively. So, inserting (\ref{E}), (\ref{M1}) and (\ref{M2}) into (\ref{rewNKL}), we obtain
\begin{align*}
N(K,L;r)=\frac{1}{q}\sum_{\substack{K<k\leq (1+\delta)K\\L<\ell\leq (1+\delta)L\\k\equiv r\mod{D}\\(k,q)=1}}1+\frac{(D,q)}{q}\sum_{\substack{K<k\leq (1+\delta)K\\L<\ell\leq (1+\delta)L\\k\equiv r\mod{D}\\(\ell,q)=1\\k\ell\equiv a\mod{(D,q)}}}1+O\bigg(\frac{x^{1-\eps/2}}{q}\bigg).
\end{align*}
We plug this into (\ref{aftdd}) and deduce that
\begin{align}\label{gbla}
\sum_{\substack{n\leq x\\n\equiv a\mod{q}}}\!\!\!\!\lambda(n)=\frac{1}{q}\sum_{KL\leq x}\Bigg(\sum_{\substack{K<k\leq (1+\delta)K\\L<\ell\leq (1+\delta)L\\(k,q)=1}}\!\!\chi(k)+(D,q)\!\!\!\!\sum_{\substack{K<k\leq (1+\delta)K\\L<\ell\leq (1+\delta)L\\(\ell,q)=1\\k\ell\equiv a\mod{(D,q)}}}\!\!\!\chi(k)\Bigg)\!+O_{\eps}\bigg(\frac{Dx^{1-\eps/2}}{q}\bigg).
\end{align}

Pursuing a reasoning similar to the one developed above for the estimation of the big-Oh term in (\ref{sNs}), we can ``intertwine" the variables $k$ and $\ell$ in the sums at the right-hand side of (\ref{gbla}). This way, one can particularly infer that
\begin{align}\label{tins}
\sum_{\substack{n\leq x\\n\equiv a\mod{q}}}\lambda(n)=\frac{1}{q}\sum_{\substack{k\ell\leq x\\(k,q)=1}}\chi(k)+\frac{(D,q)}{q}\sum_{\substack{k\ell\leq x\\(\ell,q)=1\\k\ell\equiv a\mod{(D,q)}}}\chi(k)+O_{\eps}\bigg(\frac{Dx^{1-\eps/2}}{q}\bigg).
\end{align}

Let us now focus on the second sum of the right-hand side. By the hyperbola method, we have that
\begin{align}\label{hmipop}
\sum_{\substack{k\ell\leq x\\(\ell,q)=1\\k\ell\equiv a\mod{(D,q)}}}\chi(k)=\sum_{k\leq \sqrt{x}}\chi(k)\sum_{\substack{\ell\leq x/k\\(\ell,q)=1\\ \ell\equiv ak^*\mod{(D,q)}}}1+\sum_{\substack{\ell\leq \sqrt{x}\\(\ell,q)=1}}\sum_{\substack{\sqrt{x}<k\leq x/\ell\\k\equiv a\ell^*\mod{(D,q)}}}\chi(k).
\end{align}

A routine use of M\"{o}bius inversion gives
\begin{align}\label{fnp}
\sum_{\substack{\ell\leq x/k\\(\ell,q)=1\\ \ell\equiv ak'\mod{(D,q)}}}1=\sum_{d\mid q}\mu(d)\sum_{\substack{\ell'\leq \frac{x}{dk}\\d\ell'\equiv ak^*\mod{(D,q)}}}1.
\end{align}
But, since $a$ and $k^*$ are coprime to $(D,q)$, the congruence condition at the inner sum of the right-hand side of (\ref{fnp}) is satisfied if and only if $(d,(D,q))=1$. Consequently, we have that
\begin{align}\label{mocalc}
\begin{split}
&\sum_{\substack{\ell\leq x/k\\(\ell,q)=1\\ \ell\equiv ak'\mod{(D,q)}}}1=\sum_{\substack{d\mid q\\(d,(D,q))=1}}\mu(d)\sum_{\substack{\ell'\leq \frac{x}{dk}\\\ell'\equiv ad^*k^*\mod{(D,q)}}}1\\
=&\,\frac{x}{k(D,q)}\sum_{\substack{d\mid q\\(d,(D,q))=1}}\frac{\mu(d)}{d}+O(\tau(q))=\frac{x\phi(q)}{kq\phi((D,q))}+O(\tau(q)).
\end{split}
\end{align}

So,
\begin{align}\label{1stun}
\sum_{k\leq \sqrt{x}}\chi(k)\sum_{\substack{\ell\leq x/k\\(\ell,q)=1\\ \ell\equiv ak^*\mod{(D,q)}}}1=\frac{x\phi(q)}{q\phi((D,q))}\sum_{k\leq \sqrt{x}}\frac{\chi(k)}{k}+O(\tau(q)\sqrt{x}).
\end{align}

Having dealt with the first sum at the right-hand side of (\ref{hmipop}), we turn our attention to the other. For its study, we distinguish two cases. The first case is when $(D,q)=D$, or equivalently, when $D\mid q$, whereas the second one is concerned with the scenario $(D,q)\neq D$ which is equivalent to $D\nmid q$.

$\bullet\,$ \textbf{Case 1: $D\mid q$.} In this case, we observe that for $k\equiv a\ell^*\mod{(D,q)}$, we have that $\chi(k)=\chi(a)\chi(\ell^*)=\chi(a)\chi(\ell)^{-1}=\chi(a)\chi(\ell)$ because the character is real. Hence,
\begin{align*}
\sum_{\substack{\ell\leq \sqrt{x}\\(\ell,q)=1}}\sum_{\substack{\sqrt{x}<k\leq x/\ell\\k\equiv a\ell^*\mod{(D,q)}}}\chi(k)&=\chi(a)\sum_{\substack{\ell\leq \sqrt{x}\\(\ell,q)=1}}\chi(\ell)\sum_{\substack{\sqrt{x}<k\leq x/\ell\\k\equiv a\ell^*\mod{(D,q)}}}1\\
&=\frac{\chi(a)x}{(D,q)}\sum_{\substack{\ell\leq x\\(\ell,q)=1}}\frac{\chi(\ell)}{\ell}+O(\sqrt{x}).
\end{align*}

$\bullet\,$ \textbf{Case 2: $D\nmid q$.} For this case, we apply Lemma \ref{csaps} with $(D,q)$ in place of $D'$. Hence,
\begin{align*}
\sum_{\substack{\ell\leq \sqrt{x}\\(\ell,q)=1}}\sum_{\substack{\sqrt{x}<k\leq x/\ell\\k\equiv a\ell^*\mod{(D,q)}}}\chi(k)\ll \frac{D\sqrt{x}}{(D,q)}\leq D\sqrt{x}.
\end{align*}

Merging the two cases together, we find that
\begin{align}\label{2ndun}
\sum_{\substack{\ell\leq \sqrt{x}\\(\ell,q)=1}}\sum_{\substack{\sqrt{x}<k\leq x/\ell\\k\equiv a\ell^*\mod{(D,q)}}}\chi(k)=\mathds{1}_{D\mid q}\chi(a)\frac{x}{(D,q)}\sum_{\substack{\ell\leq x\\(\ell,q)=1}}\frac{\chi(\ell)}{\ell}+O(D\sqrt{x}),
\end{align}
which concludes the study of the double sums that are present at the right-hand side of (\ref{hmipop}).

Now, we insert (\ref{1stun}) and (\ref{2ndun}) into (\ref{hmipop}) and obtain
\begin{align*}
\sum_{\substack{k\ell\leq x\\(\ell,q)=1\\k\ell\equiv a\mod{(D,q)}}}\chi(k)=\frac{x\phi(q)}{q\phi((D,q))}\sum_{k\leq \sqrt{x}}\frac{\chi(k)}{k}+\mathds{1}_{D\mid q}\chi(a)\frac{x}{(D,q)}\sum_{\substack{\ell\leq x\\(\ell,q)=1}}\frac{\chi(\ell)}{\ell}+O(D\tau(q)\sqrt{x}).
\end{align*}

Note that the divisor bound and the assumption $q\leq x^{2/3-\eps}$ imply that $D\tau(q)\sqrt{x}$ is bounded by $Dx^{1-\eps/2}/q$. Therefore, substituting the last asymptotic formula into (\ref{tins}), we deduce that
\begin{align}\label{VW}
\sum_{\substack{n\leq x\\n\equiv a\mod{q}}}\lambda(n)=U+\mathds{1}_{D\mid q}\chi(a)W+O_{\eps}\bigg(\frac{Dx^{1-\eps/2}}{q}\bigg),
\end{align}
where
\begin{align*}
U\vcentcolon=&\frac{1}{q}\sum_{\substack{k\ell\leq x\\(k,q)=1}}\chi(k)+\frac{x(D,q)\phi(q)}{q^2\phi((D,q))}\sum_{k\leq \sqrt{x}}\frac{\chi(k)}{k}\quad\text{and}\\
W\vcentcolon=&\frac{x}{q}\sum_{\substack{\ell\leq x\\(\ell,q)=1}}\frac{\chi(\ell)}{\ell}.
\end{align*}
Notice that $U$ and $W$ are independent of $a$. 

We can now find compact expressions of $U$ and $W$ in terms of the function $\lambda$. Indeed, since (\ref{VW}) was proved for an arbitrary $a$ coprime to $q$, we can sum both its sides over all the reduced residue classes $a\mod{q}$. Then, since $U$ and $W$ are independent of $a$, we make use of the orthogonality of $\chi$ and infer that
\begin{align}\label{V}
U=\frac{1}{\phi(q)}\sum_{\substack{n\leq x\\(n,q)=1}}\lambda(n)+O_{\eps}\bigg(\frac{Dx^{1-\eps/2}}{q}\bigg).
\end{align}

The term $W$ is present in (\ref{VW}) only when $D\mid q$. In that case, in order to find a concise form of $W$, we multiply both sides of (\ref{VW}) by $\chi(a)$ and then sum over all the reduced residue classes $a\mod{q}$. Doing so, the orthogonality of the character implies that
\begin{align*}
W=\frac{1}{\phi(q)}\sum_{\substack{n\leq x\\(n,q)=1}}\lambda(n)\chi(n)+O_{\eps}\bigg(\frac{Dx^{1-\eps/2}}{q}\bigg).
\end{align*}
Observe that $\lambda(n)\chi(n)=\sum_{d\mid n}(\chi(d)\chi(n))=\sum_{d\mid n}\chi(d)^2\chi(n/d)=\sum_{d\mid n}\chi(n/d)=\lambda(n)$, and so
\begin{align}\label{W}
W=\frac{1}{\phi(q)}\sum_{\substack{n\leq x\\(n,q)=1}}\lambda(n)+O_{\eps}\bigg(\frac{Dx^{1-\eps/2}}{q}\bigg).
\end{align}

We are now in position to conclude the proof of Proposition \ref{slaps}. First, we plug (\ref{V}) and (\ref{W}) into (\ref{VW}) to deduce that 
$$\sum_{\substack{n\leq x\\n\equiv a\mod{q}}}\lambda(n)=\frac{1+\mathds{1}_{D\mid q}\chi(a)}{\phi(q)}\sum_{\substack{n\leq x\\(n,q)=1}}\lambda(n)+O_{\eps}\bigg(\frac{Dx^{1-\eps/2}}{q}\bigg).$$
Then, we apply Lemma \ref{csl} to the sums $\sum_{n\leq x,(n,q)=1}\lambda(n)$ at the right, and this concludes the proof of the proposition.

\section{Proof of Theorem \ref{maint} -- Part I: Estimation of the difference $\Delta$}\label{1stpot}
Now that we have established all our auxiliary lemmas and the proof of Proposition \ref{slaps} is settled, we embark on the proof of Theorem \ref{maint}. From this point onward, the parameter $z$ is the same as the one from the statement of Lemma \ref{joka}, that is,
\begin{align}\label{z}
z\vcentcolon=D^{\min\big\{\sqrt{V/(\log \eta)},\,2\big\}},
\end{align}
where $V=\log x/\log D$. In all the forthcoming applications of Lemma \ref{auxsf}, we shall also be using the parameter
\begin{align}\label{u}
u\vcentcolon=\frac{\eps V}{200\min\{\sqrt{V/\log\eta},2\}},
\end{align}
which is defined in such a way so that $z^u=x^{\eps/200}$. The exponent $\eps$ here is the same $\eps$ that we have in the statement of Theorem \ref{maint}. Note that for $\eta$ sufficiently large in terms of $\eps$, the parameter $u$ can become as large as we please. So, whenever $u$ is involved in what follows, we shall assume that it is as large as needed. This is allowed because in Theorem \ref{maint}, we basically assume that $\eta$ is large enough with respect to $\eps$.  

In virtue of (\ref{tocon}), our goal is to bound the difference $\Delta$ defined in (\ref{Delta}) and the sums $S_1$ and $S_2$ given in (\ref{S1}) and (\ref{S2}), respectively. To this end, we shall first focus on the sifted sums of $\lambda'$ over the arithmetic progression $a\mod{q}$.

Like in Section \ref{sl'sec}, in what follows, the symbol $\bar{m}$ denotes the multiplicative inverse of the integer $m\mod{q}$.

\vspace{-10mm}

\subsection{The sifted sums of $\lambda'$ over the arithmetic progression $a\mod{q}$}
In order to study the sums of $\lambda'$ over the $z$-rough integers of the arithmetic progression $a\mod{q}$, we are applying Lemma \ref{auxsf}. In particular, utilizing the lemma with $f_1=\log$, $f_2=\chi$, and $f_3$ equal to the neutral element of Dirichlet convolution, we get that
\begin{align}\label{sl'as}
\sum_{\substack{n\leq x\\n\equiv a\mod{q}\\P^-(n)>z}}\lambda'(n)=&\sum_{\substack{d_1,d_2\leq x^{\eps/200}\\(d_1d_2,q)=1}}w(d_1)w(d_2)\sum_{\substack{k\ell \leq x/(d_1d_2)\\k\ell \equiv a\bar{d_1d_2}\mod{q}}}\log(kd_1)\chi(\ell d_2)\nn
&+O\bigg(\frac{xe^{-u/15}}{\phi(q)}\bigg(\frac{\log x}{\log z}\bigg)^{\!\!12}\,\bigg)\nn
=&\sum_{\substack{d_1,d_2\leq x^{\eps/200}\\(d_1d_2,q)=1}}w(d_1)w(d_2)\chi(d_2)\!\!\sum_{\substack{n\leq x/(d_1d_2)\\n\equiv a\bar{d_1d_2}\mod{q}}}\!\!\lambda'(n)\nn
&+\sum_{\substack{d_1,d_2\leq x^{\eps/200}\\(d_1d_2,q)=1}}w(d_1)\log d_1w(d_2)\chi(d_2)\!\!\sum_{\substack{n\leq x/(d_1d_2)\\n\equiv a\bar{d_1d_2}\mod{q}}}\!\!\lambda(n)\nn
&+O\bigg(\frac{xe^{-u/15}}{\phi(q)}\bigg(\frac{\log x}{\log z}\bigg)^{\!\!12}\,\bigg).
\end{align}

We now use Lemma \ref{l's} for the sums of $\lambda'$ at the third line of (\ref{sl'as}). Since the character is real, meaning that $\chi(d_2)\chi(a\bar{d_1d_2})\!=\chi(a)\chi(d_1)$, a careful application of Lemma \ref{l's} shall in fact yield
\begin{align}\label{mess}
\begin{split}
&\sum_{\substack{d_1,d_2\leq x^{\eps/200}\\(d_1d_2,q)=1}}w(d_1)w(d_2)\chi(d_2)\!\!\sum_{\substack{n\leq x/(d_1d_2)\\n\equiv a\bar{d_1d_2}\mod{q}}}\!\!\lambda'(n)\\
&=\frac{1}{\phi(q)}\sum_{\substack{d_1,d_2\leq x^{\eps/200}\\(d_1d_2,q)=1}}w(d_1)w(d_2)\chi(d_2)\sum_{\substack{n\leq x/(d_1d_2)\\(n,q)=1}}\!\!\lambda'(n)\\
&-\frac{\mathds{1}_{D\mid q}\chi(a)}{\phi(q)}\!\!\sum_{\substack{d_1,d_2\leq x^{\eps/200}\\(d_1d_2,q)=1}}w(d_1)\chi(d_1)w(d_2)\sum_{\substack{n\leq x/(d_1d_2)\\(n,q)=1}}\!\!\lambda'(n)\\
&+\prod_{p\mid q}\bigg(1-\frac{\chi(p)}{p}\bigg)\frac{\mathds{1}_{D\mid q}\chi(a)xL(1,\chi)}{q}\sum_{\substack{d_1,d_2\leq x^{\eps/200}\\(d_1d_2,q)=1}}\frac{w(d_1)\chi(d_1)w(d_2)}{d_1d_2}\big(\log (x/d_1d_2)-1\big)\\
&+O_{\eps}\bigg(\frac{x^{1/2+2\eps}D^{1/2+\eps}}{q^{1/4}}\bigg).
\end{split}
\end{align}

Regarding the sum at the fourth line of (\ref{mess}), we can rewrite it as
\begin{align*}
&\sum_{\substack{d_1,d_2\leq x^{\eps/200}\\(d_1d_2,q)=1}}\frac{w(d_1)\chi(d_1)w(d_2)}{d_1d_2}\big(\log (x/d_1d_2)-1\big)\\
=&\,(\log x-1)\bigg(\sum_{d_1}\frac{w(d_1)\chi(d_1)\psi_0(d_1)}{d_1}\bigg)\bigg(\sum_{d_2}\frac{w(d_2)\psi_0(d_2)}{d_2}\bigg)\\
-&\bigg(\sum_{d_1}\frac{w(d_1)\chi(d_1)\psi_0(d_1)\log d_1}{d_1}\bigg)\bigg(\sum_{d_2}\frac{w(d_2)\psi_0(d_2)}{d_2}\bigg)\\
-&\bigg(\sum_{d_1}\frac{w(d_1)\chi(d_1)\psi_0(d_1)}{d_1}\bigg)\bigg(\sum_{d_2}\frac{w(d_2)\psi_0(d_2)\log d_2}{d_2}\bigg),
\end{align*}
where $\psi_0$ is the principal character $\mod{q}$.
We repeatedly make use of Lemma \ref{FL}(b) either with $j=0$ or $j=1$ to bound all the sums over $d_1$ and $d_2$ in the parentheses above and infer that
\begin{align}\label{wchi}
\sum_{\substack{d_1,d_2\leq x^{\eps/200}\\(d_1d_2,q)=1}}\frac{w(d_1)\chi(d_1)w(d_2)}{d_1d_2}\big(\log (x/d_1d_2)-1\big)\ll \log x\prod_{\substack{p\leq z\\p\nmid q}}\bigg(1-\frac{1}{p}\bigg)\bigg(1-\frac{\chi(p)}{p}\bigg).
\end{align}

However,
\begin{align}\label{1}
\begin{split}
\prod_{\substack{p\leq z\\p\nmid q}}\bigg(1-\frac{1}{p}\bigg)&=\prod_{\substack{p\leq z\\p\mid q}}\bigg(1-\frac{1}{p}\bigg)^{-1}\prod_{p\leq z}\bigg(1-\frac{1}{p}\bigg)\\
&\leq\prod_{p\mid q}\bigg(1-\frac{1}{p}\bigg)^{-1}\prod_{p\leq z}\bigg(1-\frac{1}{p}\bigg)\ll\frac{q}{\phi(q)\log z},
\end{split}
\end{align}
and 
\begin{align}\label{2}
\begin{split}
\prod_{\substack{p\leq z\\p\nmid q}}\bigg(1-\frac{\chi(p)}{p}\bigg)&=\prod_{p\leq z}\bigg(1-\frac{\chi(p)}{p}\bigg)\prod_{p\mid q}\bigg(1-\frac{\chi(p)}{p}\bigg)^{-1}\prod_{\substack{p>z\\p\mid q}}\bigg(1-\frac{\chi(p)}{p}\bigg)\\
&\ll\prod_{p\leq z}\bigg(1-\frac{\chi(p)}{p}\bigg)\prod_{p\mid q}\bigg(1-\frac{\chi(p)}{p}\bigg)^{-1}\exp\bigg(-\sum_{\substack{p\mid q\\p>z}}\frac{\chi(p)}{p}\bigg)\\
&\leq\prod_{p\leq z}\bigg(1-\frac{\chi(p)}{p}\bigg)\prod_{p\mid q}\bigg(1-\frac{\chi(p)}{p}\bigg)^{-1}\exp\bigg(\sum_{z<p\leq x}\frac{1}{p}\bigg)\\
&\ll \frac{\log x}{\log z}\prod_{p\leq z}\bigg(1-\frac{\chi(p)}{p}\bigg)\prod_{p\mid q}\bigg(1-\frac{\chi(p)}{p}\bigg)^{-1}.
\end{split}
\end{align}

We insert (\ref{1}) and (\ref{2}) into (\ref{wchi}) and obtain that
\begin{align}
\begin{split}
&\sum_{\substack{d_1,d_2\leq x^{\eps/200}\\(d_1d_2,q)=1}}\frac{w(d_1)\chi(d_1)w(d_2)}{d_1d_2}\big(\log (x/d_1d_2)-1\big)\\
\ll&\frac{q}{\phi(q)}\bigg(\frac{\log x}{\log z}\bigg)^2\prod_{p\leq z}\bigg(1-\frac{\chi(p)}{p}\bigg)\prod_{p\mid q}\bigg(1-\frac{\chi(p)}{p}\bigg)^{-1}.
\end{split}
\end{align}

Combining this estimate with Lemma \ref{fdp}, one concludes that
\begin{align}\label{ti}
\begin{split}
&\prod_{p\mid q}\bigg(1-\frac{\chi(p)}{p}\bigg)\frac{\mathds{1}_{D\mid q}\chi(a)xL(1,\chi)}{q}\sum_{\substack{d_1,d_2\leq x^{\eps/200}\\(d_1d_2,q)=1}}\frac{w(d_1)\chi(d_1)w(d_2)}{d_1d_2}\big(\log (x/d_1d_2)-1\big)\\
\ll&\,\,\frac{x}{\eta\phi(q)}\bigg(\frac{\log x}{\log z}\bigg)^{\!2}\frac{\log D}{\log z}\leq\frac{x}{\eta\phi(q)}\bigg(\frac{\log x}{\log z}\bigg)^{\!3}.
\end{split}
\end{align}

We now replace the sums of the fourth line of (\ref{mess}) with the bound that we just established. We then proceed to insert the resulting formula into (\ref{sl'as}). This way, and after observing that the sum at the third line of (\ref{mess}) is symmetric with respect to $d_1$ and $d_2$, we deduce that
\begin{align}\label{asafc}
\begin{split}
\sum_{\substack{n\leq x\\n\equiv a\mod{q}\\P^-(n)>z}}\lambda'(n)
=&\,\,\frac{1-\mathds{1}_{D\mid q}\chi(a)}{\phi(q)}\sum_{\substack{d_1,d_2\leq x^{\eps/200}\\(d_1d_2,q)=1}}w(d_1)w(d_2)\chi(d_2)\sum_{\substack{n\leq x/(d_1d_2)\\(n,q)=1}}\!\!\lambda'(n)\\
&+\sum_{\substack{d_1,d_2\leq x^{\eps/200}\\(d_1d_2,q)=1}}w(d_1)\log d_1w(d_2)\chi(d_2)\!\!\sum_{\substack{n\leq x/(d_1d_2)\\n\equiv a\bar{d_1d_2}\mod{q}}}\!\!\lambda(n)\\
&+\,\,O_{\eps}\bigg(\frac{x}{\phi(q)}\Big(\eta^{-1}+e^{-u/15}\Big)\bigg(\frac{\log x}{\log z}\bigg)^{\!\!12}+\frac{x^{1/2+2\eps}D^{1/2+\eps}}{q^{1/4}}\,\bigg).
\end{split}
\end{align}

Our next goal before moving on to the study of the sifted sums of $\lambda'\psi_0$, where $\psi_0$ is the principal Dirichlet character $\mod{q}$, is to investigate the sums of the second line in (\ref{asafc}).

Since $x/(d_1d_2)\geq x^{1-\eps/100}$ and $\eps<1/100$, we can verify that $(x/(d_1d_2))^{2/3-\eps}\!\geq x^{58/115-\eps}\!\geq q$. Therefore, we can utilize Proposition \ref{slaps} for the sums of $\lambda$ over $a\mod{q}$, and it will in particular yield
\begin{align*}
&\sum_{\substack{d_1,d_2\leq x^{\eps/200}\\(d_1d_2,q)=1}}w(d_1)\log d_1w(d_2)\chi(d_2)\!\!\sum_{\substack{n\leq x/(d_1d_2)\\n\equiv a\bar{d_1d_2}\mod{q}}}\!\!\lambda(n)\\
=&\,\,\frac{\mathds{1}_{D\mid q}\chi(a)xL(1,\chi)}{q}\prod_{p\mid q}\bigg(1-\frac{\chi(p)}{p}\bigg)\!\!\sum_{\substack{d_1,d_2\leq x^{\frac{\eps}{200}}}}\frac{w(d_1)\chi(d_1)\psi_0(d_1)\log d_1w(d_2)\psi_0(d_2)}{d_1d_2}\\
+&\,\,\frac{xL(1,\chi)}{q}\prod_{p\mid q}\bigg(1-\frac{\chi(p)}{p}\bigg)\!\!\sum_{\substack{d_1,d_2\leq x^{\frac{\eps}{200}}}}\frac{w(d_1)\psi_0(d_1)\log d_1w(d_2)\chi(d_2)\psi_0(d_2)}{d_1d_2}+O_{\eps}\bigg(\frac{Dx^{1-\eps/3}}{q}\bigg).
\end{align*}

We are now going to apply Lemma \ref{FL}(b) four times; twice for the sums over $d_1$ and $d_2$ on the second line, and another two times for the sums over $d_1$ and $d_2$ at the third line. Combining this multiple use of Lemma \ref{FL}(b) with the estimates (\ref{1}) and (\ref{2}), as well as with Lemma \ref{fdp}, we deduce that
\begin{align*}
&\sum_{\substack{d_1,d_2\leq x^{\eps/200}\\(d_1d_2,q)=1}}w(d_1)\log d_1w(d_2)\chi(d_2)\!\!\sum_{\substack{n\leq x/(d_1d_2)\\n\equiv a\bar{d_1d_2}\mod{q}}}\!\!\lambda(n)\\
\ll_{\eps}&\,\,\frac{x(\log z)L(1,\chi)}{q}\prod_{p\mid q}\bigg(1-\frac{\chi(p)}{p}\bigg)\prod_{\substack{p\leq z\\p\nmid q}}\bigg(1-\frac{1}{p}\bigg)\bigg(1-\frac{\chi(p)}{p}\bigg)+\frac{Dx^{1-\eps/3}}{q}\\
\ll&\,\,\frac{x\log x}{\phi(q)\log z}L(1,\chi)\prod_{p\leq z}\bigg(1-\frac{\chi(p)}{p}\bigg)+\frac{Dx^{1-\eps/3}}{q}\ll \frac{x}{\eta\phi(q)}\bigg(\frac{\log x}{\log z}\bigg)^{\!2}+\frac{Dx^{1-\eps/3}}{q}.
\end{align*}

Inserting this into (\ref{asafc}), it follows that
\begin{align}\label{tfD1}
\begin{split}
\sum_{\substack{n\leq x\\n\equiv a\mod{q}\\P^-(n)>z}}\lambda'(n)
=&\,\,\frac{1-\mathds{1}_{D\mid q}\chi(a)}{\phi(q)}\sum_{\substack{d_1,d_2\leq x^{\eps/200}\\(d_1d_2,q)=1}}w(d_1)w(d_2)\chi(d_2)\sum_{\substack{n\leq x/(d_1d_2)\\(n,q)=1}}\!\!\lambda'(n)\\
&+\,\,O_{\eps}\bigg(\frac{x}{\phi(q)}\Big(\eta^{-1}+e^{-u/15}\Big)\bigg(\frac{\log x}{\log z}\bigg)^{\!\!12}\,\bigg),
\end{split}
\end{align}
where we omitted the contributions of the terms $Dx^{1-\eps/3}q^{-1}$ and $x^{1/2+2\eps}D^{1/2+\eps}q^{-1/4}$ in the big-Oh. Let us explain why. First of all, since $\eps<1/100$ and $q\leq x^{58/115-\eps}$, one can readily show that the latter expression is absorbed by the former. Then both can be ignored because $Dx^{1-\eps/3}q^{-1}$ is dominated by $x/(\eta\phi(q))$. Indeed, due to the assumption $D\leq x^{\eps/200}$ and the inequality $\eta\ll D^{50}$ which comes from Siegel's theorem, namely (\ref{siegel}), we have that $Dx^{1-\eps/3}q^{-1}\leq x^{1-\eps/4}\phi(q)^{-1}\leq xD^{-50}\phi(q)^{-1}\ll x/(\eta\phi(q))$. 

We have just completed our treatment on the sifted sums of the function $\lambda'$ over the arithmetic progression $a\mod{q}$ and we now shift our focus on the sifted sums of $\lambda'$ over the positive integers that are coprime to $q$.

\vspace{-5mm}

\subsection{The sums of $\lambda'$ over the $z$-rough integers that are coprime to $q$}

In this subsection, we aim to analyze the expression
\begin{align*}
\frac{1-\mathds{1}_{D\mid q}\chi(a)}{\phi(q)}\sum_{\substack{n\leq x\\(n,q)=1\\P^-(n)>z}}\lambda'(n),    
\end{align*}
which is subtracted in the definition of $\Delta$ in (\ref{Delta}).

Recalling that $\psi_0$ denotes the trivial character $\mod{q}$, we apply Lemma \ref{auxsf} with $a=q=1$ (these are the ones in the statement of the lemma; not to be confused with the ones in the statement of Theorem \ref{maint}), $f_1=\log\cdot\,\psi_0,f_2=\psi_0$ and $f_3$ equal to the neutral element of Dirichlet convolution. Then
\begin{align}\label{sl'c}
\sum_{\substack{n\leq x\\(n,q)=1\\P^-(n)>z}}\lambda'(n)&=\sum_{\substack{d_1,d_2\leq x^{\eps/200}\\(d_1d_2,q)=1}}w(d_1)w(d_2)\sum_{\substack{k\ell \leq x/(d_1d_2)\\(k\ell,q)=1}}\log(kd_1)\chi(\ell d_2)+O\bigg(xe^{-u/15}\bigg(\frac{\log x}{\log z}\bigg)^{\!12}\bigg)\nn
&=\sum_{\substack{d_1,d_2\leq x^{\eps/200}\\(d_1d_2,q)=1}}w(d_1)w(d_2)\chi(d_2)\sum_{\substack{n\leq x/(d_1d_2)\\(n,q)=1}}\lambda'(n)\\
&+\sum_{\substack{d_1,d_2\leq x^{\eps/200}\\(d_1d_2,q)=1}}w(d_1)\log d_1w(d_2)\chi(d_2)\sum_{\substack{n\leq x/(d_1d_2)\\(n,q)=1}}\lambda(n)+O\bigg(xe^{-u/15}\bigg(\frac{\log x}{\log z}\bigg)^{\!12}\bigg)\nonumber.
\end{align}

We can estimate the sums of the bottom line in a similar fashion as the sums
\begin{align*}
\sum_{\substack{d_1,d_2\leq x^{\eps/200}\\(d_1d_2,q)=1}}w(d_1)\log d_1w(d_2)\chi(d_2)\!\!\sum_{\substack{n\leq x/(d_1d_2)\\n\equiv a\bar{d_1d_2}\mod{q}}}\!\!\lambda(n)
\end{align*}
from (\ref{sl'as}). The only major difference is that instead of using Proposition \ref{slaps}, we now make use of the asymptotics provided by Lemma \ref{csl}. Since the remaining core details of the argument remain the same, we are not going to repeat them here for the sake of space, but we shall obtain the bound
\begin{align}\label{3}
\sum_{\substack{d_1,d_2\leq x^{\eps/200}\\(d_1d_2,q)=1}}w(d_1)\log d_1w(d_2)\chi(d_2)\sum_{\substack{n\leq x/(d_1d_2)\\(n,q)=1}}\lambda(n)
\ll_{\eps}\frac{x}{\eta}\bigg(\frac{\log x}{\log z}\bigg)^{\!2}.
\end{align}

We put (\ref{sl'c}) together with (\ref{3}) and get
\begin{align}\label{tfD2}
\begin{split}
\frac{1-\mathds{1}_{D\mid q}\chi(a)}{\phi(q)}\sum_{\substack{n\leq x\\(n,q)=1\\P^-(n)>z}}\lambda'(n)=&\,\,\frac{1-\mathds{1}_{D\mid q}\chi(a)}{\phi(q)}\sum_{\substack{d_1,d_2\leq x^{\eps/200}\\(d_1d_2,q)=1}}w(d_1)w(d_2)\chi(d_2)\sum_{\substack{n\leq x/(d_1d_2)\\(n,q)=1}}\!\!\lambda'(n)\\
&+\,\,O_{\eps}\bigg(\frac{x}{\phi(q)}\Big(\eta^{-1}+e^{-u/15}\Big)\bigg(\frac{\log x}{\log z}\bigg)^{\!\!12}\,\bigg).
\end{split}
\end{align}

This concludes the present subsection. We now have what we need to complete the estimation of the difference $\Delta$.

\subsection{Bounding $\Delta$ -- The final step}
Comparing (\ref{tfD1}) and (\ref{tfD2}) and recalling the choices of $z$ and $u$ from (\ref{z}) and (\ref{u}), respectively, we conclude that
\begin{align}\label{Delb}
\begin{split}
\Delta&=\sum_{\substack{n\leq x\\n\equiv a\mod{q}\\P^-(n)>z}}\lambda'(n)-\frac{1-\mathds{1}_{D\mid q}\chi(a)}{\phi(q)}\sum_{\substack{n\leq x\\(n,q)=1\\P^-(n)>z}}\lambda'(n)\\
&\ll_{\eps}\frac{x}{\phi(q)}\Big(\eta^{-1}+e^{-u/15}\Big)\bigg(\frac{\log x}{\log z}\bigg)^{\!\!12}\ll_{\eps}\frac{x}{\phi(q)}\bigg(\frac{u^{12}}{\eta}+e^{-u/20}\bigg)\\
&\ll_{\eps}\frac{x}{\phi(q)}\bigg(\frac{V^{12}}{\eta}+\exp\Big(-c_{\eps}\sqrt{V\log\eta}\Big)\bigg),
\end{split}
\end{align}
where $c_{\eps}$ is a positive constant depending on $\eps$.

We just finished the estimation of $\Delta$ and this completes the first big step in the proof of Theorem \ref{maint}. Our next goals are the estimations of the sums $S_1$ and $S_2$ defined in (\ref{S1}) and (\ref{S2}).

\section{Proof of Theorem \ref{maint} -- Part II: Estimation of $S_1$ and $S_2$ and completion of the proof}\label{lpot}
Let us rewrite here the definitions of $S_1$ and $S_2$.
\begin{align*}
S_1=&\sum_{\substack{k\ell\leq x,\,k>z\\(k\ell,q)=1\\P^-(k\ell)>z}}\lambda(k)\Lambda(\ell),\\
S_2=&\sum_{\substack{k\ell\leq x,\,k>z\\k\ell\equiv a\mod{q}\\P^-(k\ell)>z}}\lambda(k)\Lambda(\ell).
\end{align*}

\subsection{Bounding $S_1$}
We first bound $S_1$. By the definition of the von Mangoldt function, we have that $\Lambda(1)=0$ and that $\Lambda(n)\leq \log n$ for all $n\in\N$. Moreover, the arithmetic function $\lambda$ is non-negative and the bound $\sum_{m\leq t,P^-(m)>z}1\ll t(\log z)^{-1}$ holds for all $t>z$ (this can be derived with a direct application of Shiu's theorem \cite[Theorem 1]{sh} to the function $\mathds{1}_{P^-(\cdot)>z}$). Combining all these with Lemma \ref{joka}, we deduce that
\begin{align}\label{S1es}
\begin{split}
S_1&\leq\sum_{\substack{z<k\leq x\\P^-(k)>z}}\lambda(k)\sum_{\substack{z<\ell\leq x/k\\P^-(\ell)>z}}\Lambda(\ell)\\
&\leq \log x\sum_{\substack{z<k<x/z\\P^-(k)>z}}\lambda(k)\sum_{\substack{\ell\leq x/k\\P^-(\ell)>z}}1\ll \frac{x\log x}{\log z}\sum_{\substack{z<k<x/z\\P^-(k)>z}}\frac{\lambda(k)}{k}\\
&\ll\frac{xV}{\min\big\{\sqrt{V/\log\eta},2\big\}}\bigg(\frac{V^3}{\eta}+\exp\Big(\!\!-c\sqrt{V\log \eta}\Big)\bigg)\\
&\ll\frac{xV^4}{\eta}+x\exp\Big(-c'\sqrt{V\log\eta}\Big),
\end{split}
\end{align}
where $c$ is the constant from Lemma \ref{joka} and $c'$ is a positive absolute constant smaller than $c$. 

\subsection{Bounding $S_2$}
The treatment of $S_1$ is complete and now is the turn of $S_2$. Since $\Lambda(1)=0$ and $\Lambda(n)\leq \log n$ for all $n\in\N$, we have that
\begin{align}\label{S2ss}
S_2\leq \log x\sum_{\substack{k\ell\leq x\\k,\ell>1\\k\ell\equiv a\mod{q}\\P^-(k\ell)>z}}\lambda(k).
\end{align}

The sum of the right-hand side can be rewritten as
\begin{align}\label{arewr}
\sum_{\substack{k\ell\leq x\\k,\ell>1\\k\ell\equiv a\mod{q}\\P^-(k\ell)>z}}\lambda(k)=\sum_{\substack{k\ell\leq x\\k\ell\equiv a\mod{q}\\P^-(k\ell)>z}}\lambda(k)-\sum_{\substack{k\leq x\\k\equiv a\mod{q}\\P^-(k)>z}}\lambda(k)-\sum_{\substack{\ell\leq x\\\ell\equiv a\mod{q}\\P^-(\ell)>z}}1+1.
\end{align}

Regarding the rightmost sum in (\ref{arewr}), applying Lemma \ref{auxsf} with $f_1\equiv 1$ and $f_2$ and $f_3$ equal to the neutral element of Dirichlet convolution, we deduce that 
\begin{align}\label{bsieve}
\begin{split}
\sum_{\substack{\ell\leq x\\\ell\equiv a\mod{q}\\P^-(\ell)>z}}1&=\sum_{\substack{d\leq x^{\eps/100}\\(d,q)=1}}w(d)\sum_{\substack{\ell'\leq x/d\\\ell'\equiv a\bar{d}\mod{q}}}1+O\bigg(\frac{xe^{-u/15}(\log x)^{11}}{\phi(q)(\log z)^{12}}\bigg)\\
&=\frac{x}{q}\sum_{\substack{d\leq x^{\eps/100}\\(d,q)=1}}\frac{w(d)}{d}+O\bigg(x^{\eps/200}+\frac{xe^{-u/15}(\log x)^{11}}{\phi(q)(\log z)^{12}}\bigg).
\end{split}
\end{align}

Concerning the sums of $\lambda$ at the middle of the right-hand side of (\ref{arewr}), another application of Lemma \ref{auxsf} yields
\begin{align*}
\sum_{\substack{k\leq x\\k\equiv a\mod{q}\\P^-(k)>z}}\lambda(k)&=\sum_{\substack{d_1,d_2\leq x^{\eps/200}\\(d_1d_2,q)=1}}w(d_1)w(d_2)\sum_{\substack{mn\leq x/(d_1d_2)\\mn\equiv a\bar{d_1d_2}\mod{q}}}\chi(d_1m)+O\bigg(\frac{xe^{-u/15}(\log x)^{11}}{\phi(q)(\log z)^{12}}\bigg)\\
&=\sum_{\substack{d_1,d_2\leq x^{\eps/200}\\(d_1d_2,q)=1}}w(d_1)\chi(d_1)w(d_2)\sum_{\substack{k\leq x/(d_1d_2)\\k\equiv a\bar{d_1d_2}\mod{q}}}\lambda(k)+O\bigg(\frac{xe^{-u/15}(\log x)^{11}}{\phi(q)(\log z)^{12}}\bigg).
\end{align*}

We can now combine Proposition \ref{slaps} with Lemma \ref{FL}(b), (\ref{1}), (\ref{2}) and Lemma \ref{fdp} in order to estimate the sums at the bottom line the same way we bounded the sums
\begin{align*}
\sum_{\substack{d_1,d_2\leq x^{\eps/200}\\(d_1d_2,q)=1}}w(d_1)\log d_1w(d_2)\chi(d_2)\sum_{\substack{n\leq x/(d_1d_2)\\n\equiv a\bar{d_1d_2}\mod{q}}}\lambda(n).    
\end{align*}
The exact procedure shall eventually produce the estimate
\begin{align}\label{bstb}
\sum_{\substack{d_1,d_2\leq x^{\eps/200}\\(d_1d_2,q)=1}}w(d_1)\chi(d_1)w(d_2)\sum_{\substack{k\leq x/(d_1d_2)\\k\equiv a\bar{d_1d_2}\mod{q}}}\lambda(k)\ll_{\eps}\frac{x\log x}{\eta\phi(q)(\log z)^2}.
\end{align}

The remaining part of (\ref{arewr}) that we have not treated yet is the leftmost sum of its right-hand side. For those sums, we utilize Lemma \ref{auxsf} with $f_1=\chi$ and $f_2=f_3\equiv 1$. Then
\begin{align}\label{latter}
\sum_{\substack{k\ell\leq x\\k\ell\equiv a\mod{q}\\P^-(k\ell)>z}}\!\!\lambda(k)&=\!
\sum_{\substack{d_1,d_2,d_3\leq x^{\eps/200}\\(d_1d_2d_3,q)=1}}\!\!w(d_1)w(d_2)w(d_3)\!\!\sum_{\substack{mn\ell\leq x/(d_1d_2d_3)\\mn\ell\equiv a\bar{d_1d_2d_3}\mod{q}}}\!\!\chi(d_1m)+O\bigg(\frac{xe^{-u/15}(\log x)^{11}}{\phi(q)(\log z)^{12}}\bigg)\nn
&=\!\sum_{\substack{d_1,d_2,d_3\leq x^{\eps/200}\\(d_1d_2d_3,q)=1}}\!\!\!w(d_1)\chi(d_1)w(d_2)w(d_3)\!\!\!\sum_{\substack{k\ell\leq x/(d_1d_2d_3)\\k\ell\equiv a\bar{d_1d_2d_3}\mod{q}}}\!\!\!\lambda(k)+O\bigg(\frac{xe^{-u/15}(\log x)^{11}}{\phi(q)(\log z)^{12}}\bigg).
\end{align}

We insert (\ref{latter}) along with (\ref{bsieve}) and (\ref{bstb}) into (\ref{arewr}) and infer that
\begin{align}\label{at1}
\begin{split}
\sum_{\substack{k\ell\leq x\\k,\ell>1\\k\ell\equiv a\mod{q}\\P^-(k\ell)>z}}\lambda(k)=&\sum_{\substack{d_1,d_2,d_3\leq x^{\eps/200}\\(d_1d_2d_3,q)=1}}w(d_1)\chi(d_1)w(d_2)w(d_3)\sum_{\substack{k\ell\leq x/(d_1d_2d_3)\\k\ell\equiv a\bar{d_1d_2d_3}\mod{q}}}\lambda(k)\\
&-\frac{x}{q}\sum_{\substack{d\leq x^{\eps/200}\\(d,q)=1}}\frac{w(d)}{d}+O_{\eps}\bigg(\frac{x}{\phi(q)}\Big(\eta^{-1}+e^{-u/15}\Big)\frac{(\log x)^{11}}{(\log z)^{12}}\bigg).
\end{split}
\end{align}
Notice that the term $x^{\eps/200}$ was not included in the big-Oh term because its contribution is dominated by $x/(\eta\phi(q)\log x)$. Indeed, by the assumptions $D\leq x^{\eps/200}$ and $q\leq x^{58/115-\eps}$, as well as the inequality $\eta\ll D$, which is one particular case of Siegel's theorem (\ref{siegel}), we do conclude that $x^{\eps/200}\ll x/(D\phi(q)\log x)\ll x/(\eta\phi(q)\log x)$.

One can now pursue the arguments developed from (\ref{arewr}) up to (\ref{at1}) with the congruence conditions being replaced with the coprimality conditions $(k,q)=1$ and $(\ell,q)=1$. The only significant differences that will occur reside in the applications of Lemma \ref{auxsf}, where this time, we appeal to it with $a=q=1$ (we are thus not applying the lemma with the $a$ and $q$ from the statement of Theorem \ref{maint}) and the functions $\mathds{1}_{(\cdot,q)=1}$ and $\chi\cdot\mathds{1}_{(\cdot,q)=1}$ in place of the constant function $1$ and the character $\chi$. We shall then derive that
\begin{align}\label{at2}
\begin{split}
\frac{1}{\phi(q)}\sum_{\substack{k\ell\leq x\\k,\ell>1\\(k\ell,q)=1\\P^-(k\ell)>z}}\lambda(k)=&\sum_{\substack{d_1,d_2,d_3\leq x^{\eps/200}\\(d_1d_2d_3,q)=1}}w(d_1)\chi(d_1)w(d_2)w(d_3)\frac{1}{\phi(q)}\!\sum_{\substack{k\ell\leq x/(d_1d_2d_3)\\(k\ell,q)=1}}\lambda(k)\\
&-\frac{x}{q}\sum_{\substack{d\leq x^{\eps/200}\\(d,q)=1}}\frac{w(d)}{d}+O_{\eps}\bigg(\frac{x}{\phi(q)}\Big(\eta^{-1}+e^{-u/15}\Big)\frac{(\log x)^{11}}{(\log z)^{12}}\bigg).
\end{split}
\end{align}

We compare (\ref{at1}) to (\ref{at2}) and deduce that
\begin{align}\label{at3}
\sum_{\substack{k\ell\leq x\\k,\ell>1\\k\ell\equiv a\mod{q}\\P^-(k\ell)>z}}\lambda(k)=\frac{1}{\phi(q)}\sum_{\substack{k\ell\leq x\\k,\ell>1\\(k\ell,q)=1\\P^-(k\ell)>z}}\lambda(k)+\Delta'+O_{\eps}\bigg(\frac{x}{\phi(q)}\Big(\eta^{-1}+e^{-u/15}\Big)\frac{(\log x)^{11}}{(\log z)^{12}}\bigg),
\end{align}
where
\begin{align*}
\Delta'\vcentcolon=\sum_{\substack{d_1,d_2,d_3\leq x^{\eps/200}\\(d_1d_2d_3,q)=1}}w(d_1)\chi(d_1)w(d_2)w(d_3)\bigg(\sum_{\substack{k\ell\leq x/(d_1d_2d_3)\\k\ell\equiv a\bar{d_1d_2d_3}\mod{q}}}\lambda(k)-\frac{1}{\phi(q)}\!\!\sum_{\substack{k\ell\leq x/(d_1d_2d_3)\\(k\ell,q)=1}}\lambda(k)\bigg).
\end{align*}

Since $x/(d_1d_2d_3)\geq x^{1-3\eps/200}$ for the summation variables $d_1,d_2$ and $d_3$ in $\Delta'$, we can check that $q\leq (x/(d_1d_2d_3))^{58/115-\eps/2}$, which guarantees the applicability of Lemma \ref{compa} with $\eps/2$ in place of $\eps$ to all the differences in the sum of $\Delta'$. We thus conclude that
\begin{align}\label{fat1}
\Delta'\ll_{\eps}\frac{Dx^{1-\eps/12+3\eps/200}}{q}\leq\frac{x^{1-\eps/12+\eps/50}}{q}<\frac{x^{1-\eps/20}}{q}\leq\frac{x}{D^{10}\phi(q)}\ll \frac{x}{\eta\phi(q)},
\end{align}
where the last step follows from Siegel's theorem, that is, the bound (\ref{siegel}). 

In addition to (\ref{fat1}), the estimation of the sums from the right-hand side of (\ref{at3}) has been carried out in (\ref{S1es}), and specifically, from the second line of (\ref{S1es}) and below. Hence,
\begin{align}\label{fat2}
\frac{1}{\phi(q)}\sum_{\substack{k\ell\leq x\\k,\ell>1\\(k\ell,q)=1\\P^-(k\ell)>z}}\lambda(k)\ll\frac{x}{\phi(q)}\bigg(\frac{V^4}{\eta}+\exp\Big(-c'\sqrt{V\log\eta}\Big)\bigg),
\end{align}
where $c'$ is the same positive absolute constant that we have in (\ref{S1es}).

We now plug (\ref{fat1}) and (\ref{fat2}) into (\ref{at3}) and obtain that
\begin{align*}
\sum_{\substack{k\ell\leq x\\k,\ell>1\\k\ell\equiv a\mod{q}\\P^-(k\ell)>z}}\lambda(k)\ll_{\eps}\frac{x}{\phi(q)}\Big(V^4\eta^{-1}+e^{-u/15}+\exp\big(-c'\sqrt{V\log\eta}\big)\Big)\frac{(\log x)^{11}}{(\log z)^{12}},
\end{align*}
which is in turn put into (\ref{S2ss}) and yields the bound
\begin{align}\label{S2tfb}
\begin{split}
S_2&\ll_{\eps} \frac{x}{\phi(q)}\Big(V^4\eta^{-1}+e^{-u/15}+\exp\big(-c'\sqrt{V\log\eta}\big)\Big)\bigg(\frac{\log x}{\log z}\bigg)^{\!12}\\
&\ll_{\eps}\frac{x}{\phi(q)}\bigg(\frac{V^4u^{12}}{\eta}+e^{-u/20}+u^{12}\exp\Big(-c'\sqrt{V\log\eta}\Big)\bigg)\\
&\ll_{\eps}\frac{x}{\phi(q)}\bigg(\frac{V^{16}}{\eta}+\exp\Big(-c_{\eps}'\sqrt{V\log\eta}\Big)\bigg),    
\end{split}
\end{align}
where $c_{\eps}'>0$ is a constant that depends on $\eps$.

\subsection{Putting the pieces together -- Finalization of the proof}
Now that the estimation of both $S_1$ and $S_2$ is complete, we are in position to conclude the proof of Theorem \ref{maint}. In particular, upon consulting the exact definitions of $\Delta, S_1$ and $S_2$ from (\ref{Delta}), (\ref{S1}) and (\ref{S2}), respectively, a simple substitution of the bounds (\ref{Delb}), (\ref{S1es}) and (\ref{S2tfb}) into (\ref{tocon}) implies that
\begin{align*}
\psi(x;q,a)-\frac{\psi(x)}{\phi(q)}(1-\mathds{1}_{D\mid q}\chi(a))\ll_{\eps}\frac{x}{\phi(q)}\bigg(\frac{V^{16}}{\eta}+\exp\Big(-C_{\eps}\sqrt{V\log\eta}\Big)\bigg),
\end{align*}
where $C_{\eps}$ is a positive constant depending on $\eps$. Combining the estimate that we just obtained with the classical bound of Chebyshev $\psi(x)\gg x$, we complete the proof of Theorem \ref{maint}.



\begin{thebibliography}{99}
\bibitem{fin} J. B. Friedlander and H. Iwaniec,\ {\it A note on Dirichlet $L$-functions.} Expo. Math. 36,\, no.3--4 (2018): 343--350.
\bibitem{fip} J. B. Friedlander and H. Iwaniec,\, {\it Exceptional characters and prime numbers in arithmetic progressions.} Int. Math. Res. Not. 2003,\, no. 37,\, 2033--2050.
\bibitem{fip2} J. B. Friedlander and H. Iwaniec,\, {\it Incomplete Kloosterman sums and a divisor problem.} With an appendix by Bryan J. Birch and Enrico Bombieri, Ann. of Math. (2) 121 (1985),\, no. 2,\, 319--350.
Publications. American Mathematical Society,\, Providence,\, RI,\, 2010.
\bibitem{heath} D. R. Heath-Brown,\, {\it Prime twins and Siegel zeros.} Proc. London Math. Soc. (3), 47 (1983),\, no. 2,\, 193--224.
\bibitem{hbd3} D. R. Heath-Brown,\, {\it The divisor function $d_3(n)$ in arithmetic progressions.} Acta Arith. 47 (1986),\, no. 1,\, 29--56.
\bibitem{iwko} H. Iwaniec and E. Kowalski,\, {\it Analytic number theory.} Vol 53. American Mathematical Society,\, Providence,\, 2004.
\bibitem{js} M. Jaskari and S. Sachpazis,\, {\it The Chowla conjecture and Landau-Siegel zeroes.} Math. Proc. Camb. Phil. Soc. 2025; 179(1): 167--187.
\bibitem{oldk} D. Koukoulopoulos,\, {\it On multiplicative functions which are small on average.} Geom. Funct. Anal. 23 (2013),\, no. 5,\, 1569--1630.
\bibitem{dimb} D. Koukoulopoulos,\, {\it The distribution of prime numbers.} Graduate Studies in Mathematics, 203.\, American Mathematical Society,\, Providence,\, RI,\, 2019.
\bibitem{mm} K. Matom\"{a}ki and J. Merikoski,\, {\it Siegel zeros, twin primes, Goldbach's conjecture, and primes in short intervals.} Int. Math. Res. Not. IMRN 2023,\, no. 23,\, 20337--20384.
\bibitem{mine} S. Sachpazis,\, {\it On multiplicative functions with small partial sums.} Int. Math. Res. Not. IMRN 2024,\, no. 4,\, 2937--2964.
\bibitem{sh} P. Shiu,\, {\it A Brun-Titchmarsh theorem for multiplicative functions.} J. Reine Angew. Math. 313 (1980),\, 161--170.
\bibitem{wr} T. Wright,\, {\it Prime Distribution and Siegel Zeroes.} Preprint, \url{https://arxiv.org/abs/2311.12470v2}.
\end{thebibliography}
\end{document}